\documentclass[11pt,twoside]{article}

\usepackage{fullpage}

\usepackage{amsfonts,graphicx}
\usepackage{amsmath,amssymb,amsthm}
\usepackage{fullpage}
\usepackage{epsfig}
\usepackage{epstopdf}
\usepackage{fancyheadings}
\usepackage{graphics}
\usepackage{amsfonts}
\usepackage{amsmath}
\usepackage{psfrag} \usepackage{enumerate} \usepackage{setspace}

\usepackage{algorithm}
\usepackage{algorithmicx}
\usepackage{algpseudocode}
\usepackage{color}
\usepackage{tikz}
\usetikzlibrary{calc}
\usepackage{float}
\usepackage{url}

\theoremstyle{plain}

\newtheorem{theo}{Theorem}[section]

\newcounter{lm}
\newtheorem{lem}[lm]{Lemma}

\newcounter{pr}
\newtheorem{prop}[pr]{Proposition}
\newcounter{cr}
\newtheorem{cor}[cr]{Corollary}

\theoremstyle{definition} 

\newcounter{df}
\newtheorem{nota}{Notation}[section]
\newtheorem{de}[df]{Definition}
\newtheorem{exa}{Example}[section]
\newtheorem{as}{Assumption}[section]
\newcounter{al}
\newtheorem{alg}[al]{Algorithm}

\newcommand{\btheo}{\begin{theo}}
\newcommand{\bde}{\begin{de}}
\newcommand{\ble}{\begin{lem}}
\newcommand{\bpr}{\begin{prop}}
\newcommand{\bno}{\begin{nota}}
\newcommand{\bex}{\begin{exa}}
\newcommand{\bcor}{\begin{cor}}
\newcommand{\spro}{\begin{proof}}
\newcommand{\bas}{\begin{as}}
\newcommand{\balg}{\begin{alg}}

\newcommand{\etheo}{\end{theo}}
\newcommand{\ede}{\end{de}}
\newcommand{\ele}{\end{lem}}
\newcommand{\epr}{\end{prop}}
\newcommand{\eno}{\end{nota}}
\newcommand{\eex}{\end{exa}}
\newcommand{\ecor}{\end{cor}}
\newcommand{\fpro}{\end{proof}}
\newcommand{\eas}{\end{as}}
\newcommand{\ealg}{\end{alg}}

\theoremstyle{plain}

\newtheorem{theos}{Theorem}
\newtheorem{props}{Proposition}
\newtheorem{lems}{Lemma}
\newtheorem{cors}{Corollary}

\theoremstyle{definition}
\newtheorem{exas}{Example}
\newtheorem{algs}{Algorithm}
\newtheorem{asss}{Assumption}
\newtheorem{defns}{Definition}

\newcommand{\btheos}{\begin{theos}}
\newcommand{\etheos}{\end{theos}}
\newcommand{\bprops}{\begin{props}}
\newcommand{\eprops}{\end{props}}
\newcommand{\bdes}{\begin{defns}}
\newcommand{\edes}{\end{defns}}
\newcommand{\blems}{\begin{lems}}
\newcommand{\elems}{\end{lems}}
\newcommand{\bcors}{\begin{cors}}
\newcommand{\ecors}{\end{cors}}
\newcommand{\bexs}{\begin{exas}}
\newcommand{\eexs}{\end{exas}}
\newcommand{\balgs}{\begin{algs}}
\newcommand{\ealgs}{\end{algs}}
\newcommand{\bass}{\begin{asss}}
\newcommand{\eass}{\end{asss}}

\newcommand{\Opt}{\text{Opt} }

\newcommand{\convhull}{\mbox{convhull } }

\newcounter{pb}
\newtheorem{problem}[pb]{Problem}

\newcounter{rk}
\newtheorem{rem}[rk]{Remark}

\newcommand{\sphere}{\mathbb{S}}

\long\def\comment#1{}

\newcounter{pro}
\newtheorem{property}[pro]{Property}

\newcommand{\HACKPROOF}{\begin{proof}}
\newcommand{\HACKENDPROOF}{\end{proof}}

\DeclareMathOperator*{\argmin}{arg\,min}

\newcommand{\dist}{d}

\newcommand{\conv}{\ensuremath{\operatorname{conv}}}

\newcommand{\proj}{\Pi}

\newcommand{\ball}{\mathbb{B}}
\newcommand{\calA}{\mathcal{A}}
\newcommand{\R}{\mathbb{R}}
\newcommand{\calS}{\mathcal{S}}

\makeatletter
\long\def\@makecaption#1#2{
        \vskip 0.8ex
        \setbox\@tempboxa\hbox{\small {\bf #1:} #2}
        \parindent 1.5em  
        \dimen0=\hsize
        \advance\dimen0 by -3em
        \ifdim \wd\@tempboxa >\dimen0
                \hbox to \hsize{
                        \parindent 0em
                        \hfil 
                        \parbox{\dimen0}{\def\baselinestretch{0.96}\small
                                {\bf #1.} #2
                                } 
                        \hfil}
        \else \hbox to \hsize{\hfil \box\@tempboxa \hfil}
        \fi
        }
\makeatother

\begin{document}

\begin{center}
  {\LARGE{\bf{ Data Driven Stability Analysis of\\ Black-box Switched Linear Systems}}}

  \vspace{1cm}

  {\large
\begin{tabular}{ccc}
Joris Kenanian$^\star$, Ayca Balkan$^\star$, Raphael M. Jungers$^\dagger$, Paulo Tabuada$^\star$
\end{tabular}
}

  \vspace{.5cm}

  \texttt{jkenanian@engineering.ucla.edu,abalkan@ucla.edu, \\raphael.jungers@uclouvain.be, tabuada@ucla.edu} \\

  \vspace{.5cm}

  {\large $^\star$Department of Electrical and Computer Engineering at \\ University of California, Los Angeles (UCLA)} \\
\vspace{.1cm}

  {\large $^\dagger$UCLouvain - ICTEAM Institute, Belgium} \\

  \vspace{.5cm}

\today
\end{center}

\vspace*{.2cm}

\begin{abstract}
\noindent
Can we conclude the stability of an unknown dynamical system from the knowledge of a finite number of snapshots of trajectories? We tackle this black-box problem for switched linear systems. We show that, for any given random set of observations, one can give probabilistic stability guarantees. The probabilistic nature of these guarantees implies a trade-off between their quality and the desired level of confidence. We provide an explicit way of computing the best stability-like guarantee, as a function of both the number of observations and the required level of confidence. Our proof techniques rely on geometrical analysis, chance-constrained optimization, and stability analysis tools for switched systems, including the joint spectral radius.
\end{abstract}


\section{Introduction}
Most of the existing work on stability of dynamical systems is model-based, i.e., it requires the knowledge of a model for the considered system. Although natural in many contexts, a model may not always be available. Cyber-physical systems are an illustration of such difficulty: they consist of a large number of components of different nature (modeled by differential equations, difference equations, hybrid automata, lookup tables, custom switching logic, low-level legacy code, etc.) engaged in complex interactions with each other. Closed-form models for these complex and heterogeneous systems are equally complex or even not available, and therefore one cannot use model-based techniques in these situations. The emphasis that industry places on simulation of such systems is then not surprising, since it is always possible to simulate them despite their complexity. This raises the question of whether one can provide formal guarantees about certain properties of these complex systems, based solely on information obtained via their simulations. We focus here on one of the most important of such properties in the context of control theory: stability.

\noindent
More formally, we consider a time-varying discrete-time dynamical system of the form
\vspace{-0.88394mm}
\begin{equation}\label{eq:dynamicalsystemGeneral}
x_{k+1} = f(k, x_k),
\end{equation}
where, $x_k \in X$ is the state of the system and $k \in \mathbb{N}$ is the time index. For the rest of the paper, we use the term \emph{black-box} to refer to systems where we do not have access to the model, i.e., to $f$, yet we can indirectly learn information about $f$ by observing traces (finite trajectories) of length $l$ (in the particular case of $l=1$, these traces (trajectories) become pairs of points $(x_k, x_{k+1})$ as defined in \eqref{eq:dynamicalsystemGeneral}). We start with the following question to serve as a stepping stone: For some $l \in \mathbb{N}_{>0}$, given $N$ traces of length $l$, $(x_{i,0},x_{i,1},\dots, x_{i,l})$, $1 \leq i \leq N$, belonging to the behavior of the system \eqref{eq:dynamicalsystemGeneral}, (i.e., $x_{i,k+1} = f(k, x_{i,k})$ for any $0 \leq k \leq l-1$ and any $1 \leq i \leq N$), what can we say about the stability of System \eqref{eq:dynamicalsystemGeneral}? 

\vspace{0.3cm}

\noindent
A potential approach to this problem is to first identify the dynamics, i.e., the function $f$, and then apply existing techniques from the model-based stability analysis literature. If System \eqref{eq:dynamicalsystemGeneral} is linear, its identification and stability analysis have been extensively studied. If $f$ is not a linear function and in particular if the system is a switched system, there are two main reasons behind our quest to directly work on system behaviors and bypass the identification phase: 
\begin{itemize}
\item Identification can potentially introduce approximation errors, and can have a high computational complexity. Again, this is the case for switched systems, for which the identification problem is NP-hard \cite{lauer};
\item Even when the function $f$ is known, in general, stability analysis is a very difficult problem \cite{stabilityHard1}.
\end{itemize}

\noindent
A fortiori, the combination of these two steps in an efficient and robust way seems far from obvious. In this work, we take a first step into more complex systems than the linear case by considering the class of switched linear systems. Although we restrict ourselves to such systems, we believe that the presented results can be extended to more general models.

\vspace{0.3cm}

\noindent
In recent years, an increasing number of researchers started addressing various verification and design problems in control of black-box systems \cite{bianchini, balkan, mitra, mitra2, kozarev2016case}. In particular, the initial idea behind this paper was influenced by the recent efforts in \cite{balkan, topcu, kapinski}, and \cite{lazar} on using simulation traces to find Lyapunov functions for systems with known dynamics. In these works, the main idea is that if one can construct a Lyapunov function candidate decreasing along several finite trajectories starting from different initial conditions, it should hopefully decrease along every other trajectory. Then, once a Lyapunov function candidate is constructed, this intuition is put to test by verifying the candidate function either via off-the-shelf tools as in \cite{topcu} and \cite{kapinski}, or via sampling-based techniques as in \cite{lazar}. This also relates to almost-Lyapunov functions introduced in \cite{liberzon}, which presents a relaxed notion of stability proved via Lyapunov functions decreasing everywhere except on a small set. These approaches cannot be directly applied to black-box systems, where we do not have access to the dynamics -as in our framework. However, they are based on the following idea that we address in this paper: By observing that a candidate Lyapunov function decreases on a large number of observations, we empirically build a certain confidence that this function is a bona-fide Lyapunov function. \emph{Can we translate this empirical observation on a finite set of points into a confidence that this Lyapunov function decreases in the whole state space?} 

\vspace{0.3cm}

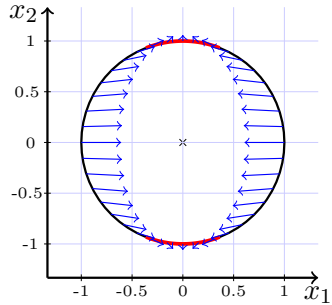
\begin{figure}
\centering
\begin{tikzpicture}[scale=0.45]
\draw[line width = 0.3mm,black,->] (-4,-4) -- (4,-4);
\draw[line width = 0.3mm,black,->] (-4,-4) -- (-4,4);

\draw (4,-4.5) node {$x_1$};
\draw (-4.7,3.8) node {$x_2$};
\draw[line width=0.1mm,black] (-0.1,-0.1) -- (0.1,0.1);
\draw[line width=0.1mm,black] (-0.1,0.1) -- (0.1,-0.1);

\foreach \i in {-2,...,2}
{
\draw[line width=0.1mm,blue!20] ({1.5*\i},-3.9) -- ({1.5*\i},4);
}

\foreach \i in {-2,...,2}
{
\draw[line width=0.1mm,blue!20] (-3.9,{1.5*\i}) -- (4,{1.5*\i});
}

\draw [line width = 0.3mm,black,domain=-68.66:68.66] plot ({3 * cos(\x)}, {3 * sin(\x)});
\draw [line width = 0.5mm,red,domain=68.66:111.34] plot ({3 * cos(\x)}, {3 * sin(\x)});
\draw [line width = 0.3mm,black,domain=111.34:248.66] plot ({3 * cos(\x)}, {3 * sin(\x)});
\draw [line width = 0.5mm,red,domain=248.66:291.34] plot ({3 * cos(\x)}, {3 * sin(\x)});

\foreach \i in {-2,...,2}
{
\draw[line width=0.1mm,black] ({1.5*\i},-4.1) -- ({1.5*\i},-3.9);
}

\draw (0,-4.4) node {\tiny{0}};
\draw (-1.5,-4.4) node {\tiny{-0.5}};
\draw (-3,-4.4) node {\tiny{-1}};
\draw (1.5,-4.4) node {\tiny{0.5}};
\draw (3,-4.4) node {\tiny{1}};

\foreach \i in {-2,...,2}
{
\draw[line width=0.1mm,black] (-4.1,{1.5*\i}) -- (-3.9,{1.5*\i});
}

\draw (-4.5,0) node {\tiny{0}};
\draw (-4.6,1.5) node {\tiny{0.5}};
\draw (-4.5,3) node {\tiny{1}};
\draw (-4.6,-1.5) node {\tiny{-0.5}};
\draw (-4.5,-3) node {\tiny{-1}};

\foreach \i in {0,...,39}
{
\draw[line width =0.15mm,->,blue] ({3*cos(\i*9)},{3*sin(\i*9)}) --  ({3*cos(\i*9)-2*0.19661*3*cos(\i*9)},{3*sin(\i*9)+2*0.03001*3*sin(\i*9)});       
}
\end{tikzpicture}
\caption{A simple dynamics and the level set of an ``almost Lyapunov function''. Even though this function decreases at almost all points in its level set, almost all trajectories diverge to infinity.}
\label{fig:levelset}
\end{figure}

\noindent
Note that, even in the case of a linear system, the connection between these two beliefs is nontrivial. In fact, one can easily construct an example where a candidate Lyapunov function decreases everywhere on its levels sets, except for an arbitrarily small subset, yet, almost all trajectories diverge to infinity. For example, the system
\[
x_{k+1} = \begin{bmatrix}
0.14 & 0\\
0 & 1.35
\end{bmatrix}x_k,
\]
admits a Lyapunov function candidate on the unit circle except on the two red areas shown in Fig. \ref{fig:levelset}. Moreover, the size of this ``violating set" can be made arbitrarily small by changing the magnitude of the unstable eigenvalue. Nevertheless, the only trajectories that do not diverge to infinity are those starting on the stable eigenspace that has zero measure.

\vspace{0.2cm}

\noindent
In this work, we take a first step for this stability inference problem, in the case of switched linear systems. In addition to the phenomenon exhibited in the above example, switched linear systems seem a priori challenging for black-box stability analysis, as both the identification and ``white-box'' stability analysis are hard for these systems. Deciding stability of a switched linear system amounts to decide whether its \emph{Joint Spectral Radius} is smaller than $1$, which is extremely hard even in the white-box setting (see, e.g., \cite{jungers_lncis}, Chapter 2, for various complexity results).

\vspace{0.2cm}

\noindent
We present an algorithm to bound the JSR of an unknown switched linear system from a finite number $N$ of observations of traces (trajectories). This algorithm partially relies on tools from the random convex optimization literature (also known as chance-constrained optimization, see \cite{campi,nemirovski,campi-garatti}), and provides an upper bound on the JSR with a user-defined confidence level. As $N$ increases, this bound gets tighter. Moreover, with a closed form expression, we characterize what is the exact trade-off between the tightness of this bound and the number of samples. In order to understand the quality of our upper bound, the algorithm also provides a deterministic lower bound. Finally, we provide a guarantee of asymptotic convergence between the upper and the lower bound, for large $N$.

\vspace{0.3cm}

\noindent
The organization of the paper is as follows: In Section~\ref{sec:preliminaries}, we introduce the problem studied and provide the necessary background in stability of switched linear systems. Then, based on finite observations for a given switched linear system, we give in Section~\ref{sec:lowerBound} a deterministic lower bound for the JSR, before presenting in Section~\ref{sec:upperbound} the main contribution of this paper, which consists in a probabilistic upper bound. We illustrate the performance of the presented techniques with some experiments in Section~\ref{sec:experiments}, and we propose future extensions of this work in Section~\ref{sec:conclusions}.

\section{Preliminaries}\label{sec:preliminaries}
\subsection{Notations}
We consider the usual Hilbert finite normed vector space $(\mathbb{R}^n,\ell_2)$, $n \in \mathbb{N}_{> 0}$, with $\ell_2$ the classical Euclidean norm. We denote by $\lVert x \rVert$ the $\ell_2$-norm of $x \in \mathbb{R}^n$. For a distance $d$ on $\mathbb{R}^n$, the distance between a set $X \subset \mathbb{R}^n$ and a point $p \in \mathbb{R}^n$ is given by $d(X,p) := \inf_{x \in X} d(x,p)$. Note that the map $p \mapsto d(X,p)$ is continuous on $\mathbb{R}^n$. Given a set $X \subset \mathbb{R}^n$, $\partial X$ denotes the boundary of set $X$.

\vspace{0.3cm}

\noindent
We also denote the set of linear functions from $\mathbb{R}^n$ to $\mathbb{R}^n$ by $\mathcal{L}(\mathbb{R}^n)$, and the set of real symmetric matrices of size $n$ by $\mathcal{S}^n$. In particular, the set of positive definite matrices is denoted by $\mathcal{S}^n_{++}$. We write $P \succ 0$ to state that $P$ is positive definite, and $P \succeq 0$ to state that $P$ is positive semi-definite. Given a set $X \subset \mathbb{R}^n$, we denote by $\wp(X)$ its powerset (i.e., the set of all its subsets), and by $X^{\mathbb{N}}$ the set of all possible sequences $(x_n)_{n \in \mathbb{N}}$, $x_n \in X$. For any $r \in \mathbb{R}_{> 0}$, we write \mbox{$rX$}$:= \{rx : x \in X\}$ to denote the scaling of ratio $r$ of $X$. We denote by $\ball$ (respectively $\sphere$) the ball (respectively sphere) of unit radius centered at the origin. We denote the ellipsoid described by the matrix $P \in \mathcal{S}^n_{++}$ as $E_P$, i.e., $E_P:= \{x \in \mathbb{R}^n: x^T P x = 1\}$. Finally, we denote the spherical projector on $\sphere$ by $\proj_{\sphere}(x) := x/\Vert x\Vert$. 

\vspace{0.3cm}

\noindent
In this paper, we only consider simple uniform probability distributions, and we believe that all the concepts can be easily intuitively understood. However, for the sake of rigor, we now develop the proper measure-theoretic setting on which our results build. We consider in this work the uniform spherical measure on $\sphere$, denoted by $\sigma^{n-1}$ ($n$ is the dimension of the space where $\sphere$ is embedded), and derived from the Lebesgue measure $\lambda$ as follows. For an ellipsoid centered at the origin, and for any of its subsets $\calA$, the \emph{sector} defined by $\calA$ is the subset $$\{t \calA, \ t \in [0,1]\} \subset\ \R^n.$$ We denote by $E_P^{\calA}$ the sector induced by $\calA \subset E_P$. In the particular case of the unit sphere, we instead write $\sphere^{\calA}$. We can notice that $E_P^{E_P}$ is the volume in $\R^n$ defined by $E_P$: $E_P^{E_P} = \{x \in \R^n: x^TPx \leq 1\}$. The spherical Borelian $\sigma$-algebra, denoted by $\mathcal{B}_{\sphere}$, is defined by $$\calA \in \mathcal{B}_{\sphere} \iff \sphere^{\calA} \in \mathcal{B}_{\mathbb{R}^n}.$$ We provide $(\sphere,\mathcal{B}_{\sphere})$ with the classical, unsigned and finite uniform spherical measure $\sigma^{n-1}$ defined by
$$\forall\ \calA \in \mathcal{B}_{\sphere},\, \sigma^{n-1}(\calA) = \frac{\lambda(\sphere^{\calA})}{\lambda(\ball)}. $$
In other words, the spherical measure of a subset of the sphere is related to the Lebesgue measure of the sector of the unit ball it induces. Notice that $\sigma^{n-1}(\sphere) = 1$.

\vspace{0.3cm}

\noindent
Since $P \in \calS_{++}^n$, we recall that it can be written in its Choleski form \begin{equation}\label{choleski}P = L^TL,\end{equation} where $L$ is an upper triangular matrix. Note that, $L^{-1}$ maps the elements of $\sphere$ to $E_P$. Then, we define the measure on the ellipsoid $\sigma_P$ on the $\sigma$-algebra $\mathcal{B}_{E_P}:=L^{-1}\mathcal{B}_\sphere$, by $\forall\, \calA \in \mathcal{B}_{E_P},\, \sigma_P({\calA}) = \sigma^{n-1}(L\calA)$. 

\vspace{0.3cm}

\noindent
For $m \in \mathbb{N}_{>0}$, we denote by $M$ the set $M = \{1,2,\dots,m\}$. The set $M$ is provided with the classical $\sigma$-algebra associated to finite sets: $\Sigma_M = \wp(M)$. We provide $(M, \Sigma_M)$ with the uniform measure $\mu_M$. For any $l \in \mathbb{N}_{>0}$, we denote by $M^l$ the $l$-Cartesian product of $M$, i.e., $M^l = \{(i_1,\dots,i_l)| i_j \in M, 1 \geq j \geq l\}$. We define $\Sigma_{M^l}$ as the product $\bigotimes^l \Sigma_M$ (which is here equal to $\wp(M)^l$), and we provide $(M^l, \Sigma_{M^l})$ with the uniform product measure $\mu_{M^l} = \otimes^l \mu_M$.

\vspace{0.3cm}

\noindent
We can now define $Z_l = \sphere \times M^l$ as the Cartesian product of $\sphere$ and $M^l$. We provide the set $Z_l$ with the product $\sigma$-algebra $\mathcal{B}_{\sphere} \bigotimes (\Sigma_{M^l})$ generated by $\mathcal{B}_{\sphere}$ and $\Sigma_{M^l}$: $\Sigma = \sigma( \pi_{\sphere}^{-1}(\mathcal{B}_{\sphere}),  \pi_{M^l}^{-1}(\Sigma_{M^l}))$, where $\pi_{\sphere}: Z_l \to \sphere$ and $\pi_{M^l}: Z_l \to M^l$ are the standard projections. On $(Z_l,\mathcal{B}_{\sphere} \bigotimes (\Sigma_{M^l}) )$, we define the product measure $\mu_l = \sigma^{n-1} \otimes \mu_{M^l}$. Note that, $\mu_l$ is the uniform probability measure on $Z_l$.

\vspace{0.3cm}

\noindent
We will also need two classical functions to compute our probabilistic upper bound, which are known as the \emph{incomplete beta function} and the \emph{regularized incomplete beta function}

\vspace{0.1cm}
\bde[\cite{handbook}, 6.6.1]
The incomplete beta function, denoted by $B$, is given by
\begin{equation*}
B: \left\{
    \begin{split}
    &\mathbb{R}_{> 0} \times \mathbb{R}_{> 0} \times \mathbb{R}_{> 0} \to \mathbb{R}_{\geq 0}\\ 
    &(x,a,b) \mapsto B(x;a,b) = \int_0^x t^{a-1} (1-t)^{b-1} dt.
    \end{split}
  \right.
\end{equation*}
\ede

\bde[\cite{handbook}, 6.6.2]
The regularized incomplete beta function, denoted by $I$, is given by
\begin{equation*}
I: \left\{
    \begin{split}
    &\mathbb{R}_{> 0} \times \mathbb{R}_{> 0} \times \mathbb{R}_{> 0} \to \mathbb{R}_{\geq 0}\\ 
    &(x,a,b) \mapsto I(x;a,b) = \frac{B(x;a,b)}{B(1;a,b)} .
    \end{split}
  \right.
\end{equation*}
\ede

\noindent
For given values of parameters $a >0$ and $b>0$, the inverse of the regularized incomplete beta function with parameters $a,b$, denoted by $I^{-1}(y;a,b)$, is the function whose output is $x>0$ such that $I(x;a,b) = y$ \cite{betafct}.

\subsection{Stability of Switched Linear Systems}\label{sec:stab}
A \emph{switched linear system}, defined by a set of modes (matrices) \mbox{$\mathcal{M}= \{A_i, i \in M \}$}, is a time-varying discrete-time dynamical system of the form \eqref{eq:dynamicalsystemGeneral}, with $f(k,x_k) = A_{\tau(k)}x_k$, that is:
\begin{equation}\label{eq:switchedSystem}
x_{k+1} = A_{\tau(k)}x_k,
\end{equation}
for any $k \in \mathbb{N}$. Here, the signal $\tau \in M^{\mathbb{N}}$ is called the \emph{switching sequence}, and can take arbitrary values in $M$. Note that such systems are homogeneous, i.e., for any $\gamma > 0$, $f(k,\gamma x_k) = \gamma f(k,x_k)$. In this paper, we assume to not have access to $\mathcal{M}$ nor to the switching sequence. The only information available is (an upper bound on) $m$, the cardinality of $\mathcal{M}$.

\noindent
We are interested in the \emph{uniform asymptotic stability} of this system, that is, we want to guarantee the following property: $$\forall \tau \in M^{\mathbb{N}}, \, \forall x_0 \in \mathbb{R}^n, \, \lVert x_k \rVert \xrightarrow[k \to \infty]{} 0.$$
The joint spectral radius of a set of matrices $\mathcal{M}$ characterizes the stability of the underlying switched linear system \eqref{eq:switchedSystem} defined by $\mathcal{M}$ \cite{jungers_lncis}. This quantity is an extension to switched linear systems of the classical spectral radius for linear systems. It is the maximum asymptotic growth rate of the norm of the state under the dynamics \eqref{eq:switchedSystem}, over all possible initial conditions and sequences of matrices of $\mathcal{M}$.

\bde[from \cite{jungers_lncis}] 
Given a finite set of matrices \mbox{$\mathcal{M} \subset \mathbb{R}^{n\times n}$}, its \emph{joint spectral radius} (JSR) is given by $$\rho(\mathcal{M}) =\lim_{k \rightarrow \infty} \max_{i_1,\dots, i_k} \left\{ ||A_{i_1} \dots A_{i_k}||^{1/k}: A_{i_j} \in \mathcal{M}\ \right\}. $$
\ede

\begin{property}[\cite{jungers_lncis}, Corollary 1.1]
Given a finite set of matrices $\mathcal{M}$, the corresponding switched dynamical system is stable if and only if $\rho(\mathcal{M})<1$.
\end{property}

\btheos[\cite{jungers_lncis}, Prop. 1.3]\label{prop:scaling}
Given a finite set of matrices $\mathcal{M}$, and any invertible matrix $T$, $\rho(\mathcal{M})=\rho(T \mathcal{M} T^{-1})$, i.e., the JSR is invariant under similarity transformations (and is a fortiori a homogeneous function: $\forall  \gamma > 0$, $\rho \left( \mathcal{M}/\gamma \right) = \rho(\mathcal{M})/\gamma$).
\etheos

\noindent
The JSR also relates to a tool classically used in control theory to study stability of systems: Lyapunov functions. We will consider here a family of such functions that is particularly adapted to the case of switched linear systems.

\bde
Consider a finite set of matrices $\mathcal{M} \subset \mathbb{R}^{n \times n}$. A \emph{common quadratic form (CQF)} for a system \eqref{eq:switchedSystem} with set of matrices $\mathcal{M}$, is a positive definite matrix $P \in \mathcal{S}_{++}^n$ such that for some $\gamma \geq 0$, 
\begin{equation}\label{eq:lyap}
\forall A \in \mathcal{M}, A^T P A \preceq \gamma^2P.
\end{equation}
\ede
\noindent
CQFs are useful because they can be computed, when they exist, with semidefinite programming (see \cite{boyd}), and they constitute a stability guarantee (when $\gamma < 1$, they are Lyapunov functions) for switched systems as we formalize next.

\btheos[\cite{jungers_lncis}, Prop. 2.8 and Thm. 2.11]\label{thm:cqlf} 
Consider a finite set of matrices $\mathcal{M}$.
\vspace{-0.2cm}
\begin{itemize} 
\item If there exist $\gamma \geq 0$ and $P \succ 0$ such that Equation \eqref{eq:lyap} holds, then $\rho(\mathcal{M}) \leq \gamma$.
\vspace{0.2cm}
\item If $\rho(\mathcal{M}) < \frac{\gamma}{\sqrt{n}},$ there exists a CQF, $P$, such that $\forall A \in \mathcal{M},\, A^T P A \preceq \gamma^2 P.$
\end{itemize}
\etheos

\noindent
For any $\gamma < 1$, this theorem provides both a Lyapunov and a \emph{converse Lyapunov result}: if there exists a CQF, then our system is stable; if there is, on the contrary, no such stability guarantee, one may conclude a lower bound on the JSR. We obtain then an approximation algorithm for the JSR. It turns out that one can still refine this technique, in order to improve the error factor $1/\sqrt{n}$, and asymptotically get rid of it. This is a well-known technique for the ``white-box'' computation of the JSR, which we summarize in the following corollary.
\begin{cor}\label{cor:approx-products}
Fix $\gamma \geq 0$. For any finite set of matrices such that $\rho(\mathcal{M}) < \frac{\gamma}{\sqrt[2l]{n}}$ with $\gamma \geq 0$, there exists a CQF for $\mathcal{M}^l:=\{\Pi_{j=1}^l A_{i_j}: A_{i_j} \in \mathcal{M}\}$, that is, a $P\succ 0$ such that:
\begin{equation} \label{eq:LMI}
\forall\ \mathbf{A} \in \mathcal{M}^l,\, \mathbf{A}^T P \mathbf{A} \preceq \gamma^{2l} P.
\end{equation}
\end{cor}
\HACKPROOF
It is easy to see from the definition of the JSR that $\rho(\mathcal{M}^l)=\rho(\mathcal{M})^l$. Thus, applying Theorem \ref{thm:cqlf} to the finite set $\mathcal{M}^l,$ one directly obtains the corollary.
\HACKENDPROOF

\noindent
Note that, the smaller $\gamma$ is in Theorem~\ref{thm:cqlf}, the tighter is the upper bound we get on $\rho(\mathcal{M})$. In order to properly analyze our setting, where matrices are unknown, let us reformulate \eqref{eq:LMI} in another form. For any $l \in \mathbb{N}_{>0}$, we can consider the optimal solution $\gamma^*$ of the following optimization problem:
\begin{equation}\label{eqn:campiOpt1}
\begin{aligned}
\hspace{-0.3cm}\text{min}_{\gamma, P} &\qquad \gamma \\
\hspace{-0.3cm}\text{s.t.} \quad  &(\mathbf{A} x)^T P \mathbf{A} x \leq \gamma^{2l} x^T P x, \mathbf{A} \in \mathcal{M}^l, \,\forall x \in \mathbb{R}^n\\
&P \succ 0.
\end{aligned}
\end{equation}

\noindent
Notice that we can restrict the set of constraints by restricting $x$ to $\sphere$, due to the homogeneity of the system. Homogeneity indeed implies that it is sufficient to show the decrease of a CQF on an arbitrary set enclosing the origin. Hence the optimization problem \eqref{eqn:campiOpt1} is equivalent to the optimization problem:

\begin{equation}\label{eqn:campiOpt2}
\begin{aligned}
\hspace{-0.3cm}\text{min}_{\gamma, P} &\qquad \gamma \\
\hspace{-0.3cm}\text{s.t.} \quad  &(\mathbf{A} x)^T P \mathbf{A} x \leq \gamma^{2l} x^T P x, \mathbf{A} \in \mathcal{M}^l, \,\forall x \in \sphere\\
&P \succ 0.
\end{aligned}
\end{equation}

\noindent
The above equation will provide a clear algebraic formalization of our black-box problem: our goal amounts to find a solution to a convex problem with an infinite number of constraints, while only sampling a finite number of them.

\subsection{Problem Formulation}
Let us now formally present the problem addressed in this paper. We recall that we only observe $N$ finite traces of length $l \in \mathbb{N}_{>0}$, i.e., $N$ sequences of states $(x_k,x_{k+1},\dots,x_{k+l})$ where $x_{k+i+1}$ and $x_{k+i}$ are related by \eqref{eq:switchedSystem}. Note that such sequences of states depend both on the initial state $x_k$ and the switching sequence $\tau(k)$ which is assumed to be unknown. In other words, we do not observe the mode or the matrices used to produce the trajectories. We do not have access to the process through which the system picks the modes. The user's knowledge is limited to the number of modes (or an upper bound on this number) and the dimension of the system. We assume that these trajectories are generated from a finite number of initial conditions $x_{i,0} \in \sphere$, $1 \leq i \leq N$ enumerating the observations, and that a random sequence of $l$ matrices is applied to each of these points. We randomly draw the initial conditions from $\mathbb{S}$, observe them and the $l$ subsequent state values produced by the system. Sampling the initial conditions from $\mathbb{S}$ is without loss of generality, since any trajectory in $\mathbb{R}^n$ can be rescaled so that $x_{i,0} \in \sphere$, by homogeneity of the system. To a given observed trajectory $(x_k,x_{k+1},\dots,x_{k+l})$, we can associate the corresponding probability event $(x_k,j_1,\dots,j_l)$ which is another $(l+1)$-tuple. Formally, with a fixed probability space $(\Omega, \mathcal{F}, \mathbb{P})$, we consider the random variables $X_0: \Omega \mapsto \mathbb{S}$ and $\theta_i: \Omega \mapsto M$, for $1 \leq i \leq l$, such that $X_0$ has uniform distribution, and the $\theta_i$ are independent and also have uniform distribution. Thus, to a given random finite set of trajectories of length $l$, we can associate an underlying uniform sample of $N$ such $(l+1)$-tuples in $Z_l = \sphere \times M^l$, denoted by
\begin{equation}\label{eq:omega}
\omega_N := \{(x_{i,0}, j_{i,1}, \dots, j_{i,l}), 1 \leq i \leq N \}  \subset Z_l
\end{equation}
In other words, a set of $N$ available observations $\{(x_{i,0},x_{i,1}, \dots, x_{i,l}), 1\leq i\leq N \}$ can be rewritten, for all $1 \leq i \leq N$ and $1 \leq k \leq l$, as $x_{i,k}= A_{j_{i,k}} \dots A_{j_{i,1}} x_{i,0}$, with the $j_{i,k}$ being unobserved variables that take their values in $M$.

\begin{rem}\label{rem:probGeneralize}
Let us motivate our assumption on the uniform drawing of the matrices. We assumed that we only have access to random observations of the state of the system, and ignore the process that generates them. In particular, we ignore the process that selects the modes at each time step, and model it as a random process. We suppose that with nonzero probability, each mode is active: the problem would indeed not be solvable otherwise, since the system would be unidentifiable with probability $1$ and would prevent to ever observe some of its behaviors. We take this distribution to be uniform since we cannot say that some modes are more likely a priori. Our results extend to the case where the distribution is not uniform as long as we have a nonzero lower bound on the probability of each mode.
\end{rem}

\noindent
In this work, we aim at understanding what type of guarantees one can obtain on the stability of System \eqref{eq:switchedSystem} (that is, on the JSR of $\mathcal{M}$) from the sample \eqref{eq:omega}. More precisely, we answer the following problem:
\begin{problem}\label{problem} 
Consider a finite set of matrices, $\mathcal{M},$ describing a switched system \eqref{eq:switchedSystem}, and suppose that one has a set of $N$ observations $(x_{i,0},x_{i,1},...,x_{i,l}), i=1,...,N$ corresponding to an event $\omega_N,$ sampled from $Z_l$ with the uniform measure $\mu_l$.
\begin{itemize}
\item For a given confidence level $\beta \in [0,1)$, provide an upper bound on $\rho(\mathcal{M})$, that is, a number $\overline{\rho(\omega_N)}$ such that $$\mu_l^N \left( \{\omega_N: \ \rho(\mathcal{M}) \leq \overline{\rho(\omega_N)} \} \right) \geq \beta.$$
\item For the same given level of confidence $\beta$, provide a lower bound $\underline{\rho(\omega_N)}$ on $\rho(\mathcal{M})$.
\end{itemize}
\end{problem}
\begin{rem}
We will see in Section~\ref{sec:lowerBound} that we can even derive a deterministic lower bound for a given (sufficiently high) number of observations. 
\end{rem}

\noindent
We will see in Theorem~\ref{thm:mainTheorem} that for any level of confidence $\beta$, it is always possible to provide an upper bound for Problem~\ref{problem} which tends to the JSR when the number of sampled points increases. In particular, for any (large enough) number of samples, it is always possible to provide such an upper bound that is finite. Thus, we obtain a Hypothesis test for the question 'Is the system stable?' with an a priori fixed probability of false positive (equal to the number $1 - \beta$), and a probability of false negative that tends to zero when the number of samples increases. 

\vspace{0.3cm}

\noindent
The key insight is to leverage the fact that conditions for the existence of a CQF for \eqref{eq:switchedSystem} can be obtained by considering a finite number of traces in $\mathbb{R}^n$ of the form $(x_k,x_{k+1}, \dots, x_{k+l})$. Developing that insight leads us to the following algorithm, that is the main result of our work and that solves Problem~\ref{problem}:
\begin{alg}[Probabilistic upper bound]
 \ 
\newline
\textbf{Input:} observations produced by a uniform random sample $\omega_N \subset Z_l$ of size $N \geq \frac{n(n+1)}{2}+1$;

\noindent
\textbf{Input:} $\beta$ desired level of certainty;

\textbf{Compute:} a candidate for the upper bound, $\gamma^{*}(\omega_N)$, solution of the convex optimization problem \eqref{eq:lowerbound};\\
(observe that \eqref{eq:lowerbound} does not require the explicit knowledge of the matrices ${\bf A_j}$)

\textbf{Compute:} $\varepsilon(\beta,\omega_N)$ the proportion of points where our inference on the upper bound may be invalid;

\textbf{Compute:} $\delta(\varepsilon) \leq 1$ a correcting factor; ($\delta \xrightarrow[N \to \infty]{} 1$)

\noindent
\textbf{Output:} $\frac{\gamma^{*}(\omega_N)}{\sqrt[2l]{n}} \leq \rho \leq \frac{\gamma^{*}(\omega_N)}{\sqrt[l]{\delta(\varepsilon)}}$;\\ 
(the right-hand side inequality is valid with probability at least $\beta$).
\end{alg}

\section{A Deterministic Lower Bound}\label{sec:lowerBound}
In Section~\ref{sec:stab}, we presented an optimization problem, \eqref{eqn:campiOpt2}, that provides a stability guarantee. Nevertheless, this problem has infinitely many constraints and observing a finite number of traces only gives us access to a restriction of it (with finitely many constraints). We consider then the following optimization problem:
\begin{equation}\label{eq:lowerbound}
\begin{aligned}
\hspace{-0.4cm}\text{min}_{\gamma, P} \quad \gamma \\
\text{s.t.}  &\qquad (\mathbf{A_j} x)^T P \mathbf{A_j} x \leq \gamma^{2l} x^T P x, \forall (x, \mathbf{j}) \in \omega_N\\
& \qquad P \succ 0,\ \gamma \geq 0. \\
\end{aligned}
\end{equation}
with optimal solution $\gamma^{*}(\omega_N)$, and where $\mathbf{A_j} :=  A_{j_l} A_{j_{l-1}} \dots A_{j_1}$ and $\mathbf{j}:=\{j_1,\dots, j_l\}$. Note that, \eqref{eq:lowerbound} can be efficiently solved by semidefinite programming and bisection on the variable $\gamma$ (see \cite{boyd}). Note also that solving this program can be done in practice only through the knowledge of the observations: even though the $A_{j_i}$ are not known, the program only requires the knowledge of $\mathbf{A_j} x,$ which is known through the observations. In this section, we provide a theorem for a deterministic lower bound based on the observations given by $\omega_N$, whose accuracy depends on the ``horizon'' $l$.

\btheos \label{thm:lowerbound}
For an arbitrary $l \in \mathbb{N}_{>0}$ and a given uniform sample $\omega_N \subset Z_l$, by considering $\gamma^*(\omega_N)$ the optimal solution of the optimization problem \eqref{eq:lowerbound}, we have $$\rho(\mathcal{M}) \geq \frac{\gamma^*(\omega_N)}{\sqrt[2l]{n}}.$$ 
\etheos
\vspace{-0.8cm}
\HACKPROOF
Let $\nu >0$. By definition of $\gamma^*(\omega_N)$, there exists no matrix $P \in \mathcal{S}^n_{++}$ such that, $\forall x \in \mathcal{S}, \, \forall \mathbf{A_j} \in \mathcal{M}^{l}$,
\begin{equation*}
(\mathbf{A_j} x)^T P \mathbf{A_j} x \leq (\gamma^*(\omega_N) -\nu)^{2l} x^T P x.
\end{equation*}
Taking the contrapositive of Corollary \ref{cor:approx-products}, this implies that $\rho(\mathcal{M}) \geq \frac{(\gamma^*(\omega_N) -\nu)^{l}}{\sqrt[2l]{n}}$. Since this is valid for any $\nu>0,$ we finally obtain the claim.
\HACKENDPROOF

\section{A Probabilistic Stability-like Guarantee}\label{sec:upperbound}

\begin{subsection}{A Partial Upper Bound}\label{sec:IntroMainThm}
In this section, we show how to compute an upper bound on $\rho$, with a user-defined confidence $\beta \in [0, 1)$. We do this by constructing an $l$-step CQF which is valid with probability at least $\beta$. The existence of an $l$-step CQF implies $\rho \leq \gamma^*$ due to Theorem \ref{thm:cqlf}. As we will see below, the quality of our bound will depend on geometrical properties of the CQF found; more precisely, the smaller the condition number of the corresponding matrix $P$, the better will be our bound. In practice, one can minimize the condition number of the solution $P$ in a second step, after computing $\gamma^*$ from \eqref{eq:lowerbound}.  However, for the sake of rigor and clarity of our proofs, we introduce a slighly different optimization problem. We consider for the rest of the discussion the following optimization problem, that we denote by $\Opt(\omega_N)$:
\begin{equation}\label{eqn:campiOpt03}
\begin{aligned}
\hspace{-0.4cm}&\min_{P} \quad \lambda_{\max}(P) \\
\hspace{-0.3cm}&\text{s.t.} \quad (\mathbf{A_j}x)^T P \mathbf{A_j} x \leq { \left( (1 +\eta)\gamma^*(\omega_N) \right)}^{2l} x^T P x, \ \forall\ (x, \mathbf{j}) \in \omega_N \\
& P \succeq I, \\
\end{aligned}
\end{equation}
with $\eta > 0$, and where $\gamma^*(\omega_N)$ is the optimal solution to the optimization problem \eqref{eq:lowerbound}. Let us analyze the relationship between $\Opt(\omega_N)$ and the optimization problem \eqref{eq:lowerbound}. Firstly, thanks to the homogeneity of system \eqref{eq:switchedSystem}, we can replace the constraint $P \succ 0$ in the initial problem with the constraint $P \succeq I$. Secondly, as discussed above, the objective function $\lambda_{\max}(P)$ (which is convex) can be added once $\gamma^*$ is computed, in order to minimize the condition number. Lastly, we introduced a ``regularization parameter'', $\eta > 0$, which ensures strict feasibility of $\Opt(\omega_N)$. As the reader will see, we will derive results valid for arbitrarily small values of $\eta$. This will then not hamper the practical accuracy of our technique, while allowing us to derive a theoretical asymptotic guarantee (i.e., for a large number of observations). We denote the optimal solution of $\Opt(\omega_N)$ by $P(\omega_N)$, and drop the explicit dependence of $P$ on $\omega_N$ when it is clear from the context.

\vspace{0.3cm}

\noindent
The intriguing question of whether the optimal solution of this sampled problem is a feasible solution to \eqref{eqn:campiOpt2} has been widely studied in the literature \cite{campi}. Under certain technical assumptions, one can bound the proportion of the constraints of the original problem \eqref{eqn:campiOpt2} that are violated by the optimal solution of $\Opt(\omega_N)$, with some probability which is a function of the sample size $N$. In the following theorem, we adapt a classical result from the random convex optimization literature to our problem.

\btheos[adapted from Theorem 3.3\footnotemark, \cite{campi}]\label{mainTheorem0}
Consider the optimization problem $\Opt(\omega_N)$ given in \eqref{eqn:campiOpt03}, where $\omega_N$ is a uniform random sample drawn from the set $Z_l$. Let $d = \frac{n(n+1)}{2}$ be the dimension of the decision variable $P$ of $\Opt(\omega_N)$ and $N \geq d+1$. Then, for all $\varepsilon \in (0,1]$ the following holds:
\begin{equation}\label{eqn:violation}
\mu_l^N\hspace{-1mm}\left\{ \omega_N \in Z_l^N: \mu_l \left( V(\omega_N) \right) \leq \varepsilon \right\}\hspace{-1mm} \geq \beta(\varepsilon, N),
\end{equation}
where $\mu_l^N$ denotes the product probability measure on $Z_l^N$, $\beta(\varepsilon, N) =  1- \sum_{j=0}^{d} \binom{N}{j}\varepsilon^j (1-\varepsilon)^{N-j}$, and $V(\omega_N)$ is the set $\{(x,\mathbf{j}) \in Z_l: (\mathbf{A_j} x)^T P(\omega_N) \mathbf{A_j} x > (\gamma_{\omega_N}^{*})^{2l} x^T P(\omega_N) x\}$, i.e., it is the subset of $Z_l$ for which the considered $\gamma^*$-contractivity is violated by the optimal solution of $\Opt(\omega_N)$.

\noindent
The quantity $\varepsilon$ can also be seen as a function of $\beta$ and $N$: $\varepsilon(\beta, N) = 1 - I(1-\beta, N-d,d+1)$ (see the proof of Theorem 3.3 in \cite{campi}).

\etheos

\footnotetext{Theorem 3.3 in \cite{campi} requires $\Opt(\omega_N)$ to satisfy the following technical assumptions:\begin{enumerate}
\item When the problem $\Opt(\omega_N)$ admits an optimal solution, this solution is unique.
\item Problem $\Opt(\omega_N)$ is nondegenerate with probability $1$.
\end{enumerate}
Here, the first assumption can be enforced if required by adding a tie-breaking rule to $\Opt(\omega_N)$ as explained in Appendix A in \cite{tiebreak}, while the second assumption can be lifted, as explained in PART 2b in \cite{campi-garatti}, thanks to the introduction of ``constraint heating''.
}

\begin{cor}\label{cor:gettingRidOfm}
Consider a set of matrices $\mathcal{M}$, a sample $\omega_N \subset Z_l$, $\gamma^{*}$ the optimal solution of \eqref{eq:lowerbound} and $P \succ 0$ the optimal solution of $\Opt(\omega_N)$. Then, with the notation of Theorem~\ref{mainTheorem0}, for any $\varepsilon>0$, with confidence $\beta(\varepsilon,N),$ one has:
\begin{equation}\label{eqn:P0}
(\mathbf{A_j} x)^T P \mathbf{A_j} x \leq (\gamma^{*})^{2l} x^T P x, \forall x \in \sphere \setminus \tilde{\sphere}, \forall \mathbf{j} \in M^l
\end{equation} 
where $\tilde{\sphere} \subset \sphere$, the projection of $V(\omega_N)$ on $\sphere$, has measure $\sigma^{n-1}(\tilde{\sphere}) \leq \varepsilon m^l $.
\end{cor}

\noindent
This result allows us to make abstraction of the probabilistic setting: by accepting a confidence level of $\beta$ smaller than one, we may assume that all the points except a small set satisfy the Lyapunov equation (5). The case where we have the equality $\sigma(\tilde{\sphere}) = \varepsilon m^l$ corresponds to the case where every point $x \in \tilde{\sphere}$ violates \eqref{eqn:P0} for exactly one value of $\mathbf{j}$.

\HACKPROOF
We know that $\Sigma_{M^l}$ is the disjoint union of its $2^{m^l}$ elements $\{\mathcal{M}^l_i, i \in \{1,2, \dots, 2^{m^l} \} \}$. Then $V(\omega_N)$ can be written as the disjoint union $V(\omega_N) = \sqcup_{1 \leq i \leq 2^{m^l}} (\mathcal{S}_i, \mathcal{M}^l_i)$ where $\mathcal{S}_i \in \Sigma_{\sphere}$. We notice that $\sphere' = \sqcup_{1 \leq i \leq 2^{m^l}} \mathcal{S}_i$, and
\begin{equation*}
\sigma^{n-1} (\tilde{\sphere}) = \sum_{1 \leq i \leq 2^{m^l}} \sigma^{n-1} (\mathcal{S}_i).
\end{equation*}
We have 
\begin{eqnarray*}
\mu_l(V(\omega_N)) &=& \mu_l \left( \sqcup_{1 \leq i \leq 2^{m^l}} (\mathcal{S}_i, \mathcal{M}^l_i) \right) = \sum_{1 \leq i \leq 2^{m^l}} \mu_l \left( \mathcal{S}_i, \mathcal{M}^l_i \right) \\
 &=& \sum_{1 \leq i \leq 2^{m^l}} \sigma^{n-1} \otimes \mu_{M^l} \left( \mathcal{S}_i, \mathcal{M}^l_i \right) \\
 &=& \sum_{1 \leq i \leq 2^{m^l}} \sigma^{n-1}(\mathcal{S}_i) \mu_{M^l} (\mathcal{M}^l_i).
\end{eqnarray*}

\noindent
Note that we have $\min_{(j_1,\dots,j_l) \in M^l} \mu_{M^l}(\{j_1,\dots,j_l\}) = \frac{1}{m^l}.$ 

\noindent
Then since $ \forall \ i$, $\mu_{M^l}(\mathcal{M}^l_i) \geq \min_{(j_1,\dots,j_l) \in M^l} \mu_{M^l}(\{j_1,\dots,j_l\}) = \frac{1}{m^l}$, we get:
\begin{equation}
\sigma^{n-1}(\tilde{\sphere}) \leq \frac{\mu_l(V(\omega_N))}{\frac{1}{m^l}} \leq m^l \varepsilon.
\end{equation}
This means that 
\begin{equation}\label{eqn:P2}
\begin{aligned}
& (A_{j_l} A_{j_{l-1}} \dots A_{j_1} x)^T P (A_{j_l} A_{j_{l-1}} \dots A_{j_1} x) \leq \gamma^{2l} x^T P x, \forall\, x \in \sphere \setminus \tilde{\sphere}, \,\forall\, (j_1,\dots,j_l) \in M^l,
\end{aligned}
\end{equation}
where $\sigma^{n-1}(\tilde{\sphere}) \leq m^l \varepsilon.$
\HACKENDPROOF

\noindent
The above results allow us to conclude, from a finite number of observations, that with probability $\beta$ (where $\beta$ goes to $1$ as $N$ goes to infinity), the required property is actually satisfied for the complete sphere $\sphere$, except on a small set of measure at most $\tilde{\varepsilon} = \varepsilon m^l$. This means that, the ellipsoid $E_P$ computed by $\Opt(\omega_N)$ is ``almost invariant''  except on a set of measure bounded by $\tilde{\varepsilon}$. This is represented in Fig. 1. for the case $n=2$, where the red points of $E_P$ are points that might violate the invariance constraint. Here, the set of red points has measure at most $\tilde{\varepsilon}$.
\begin{figure}[H]\label{fig:ellipsoid}
\begin{center}
\begin{tikzpicture}[scale=0.5]
\draw[line width=0.3mm,black,->] (-5,0) -- (5,0);
\draw[line width=0.3mm,black,->] (0,-3) -- (0,3);

\draw [dotted, line width = 0.5mm,red,domain=-30:30] plot ({4.5 * cos(\x)}, {2.5 * sin(\x)});
\draw [line width = 0.3mm,black,domain=30:55] plot ({4.5 * cos(\x)}, {2.5 * sin(\x)});
\draw [dotted, line width = 0.5mm,red,domain=55:70] plot ({4.5 * cos(\x)}, {2.5 * sin(\x)});
\draw [line width = 0.3mm,black,domain=70:95] plot ({4.5 * cos(\x)}, {2.5 * sin(\x)});
\draw [dotted, line width = 0.5mm,red,domain=95:100] plot ({4.5 * cos(\x)}, {2.5 * sin(\x)});
\draw [line width = 0.3mm,black,domain=100:150] plot ({4.5 * cos(\x)}, {2.5 * sin(\x)});
\draw [dotted, line width = 0.5mm,red,domain=150:210] plot ({4.5 * cos(\x)}, {2.5 * sin(\x)});
\draw [line width = 0.3mm,black,domain=210:235] plot ({4.5 * cos(\x)}, {2.5 * sin(\x)});
\draw [dotted, line width = 0.5mm,red,domain=235:250] plot ({4.5 * cos(\x)}, {2.5 * sin(\x)});
\draw [line width = 0.3mm,black,domain=250:275] plot ({4.5 * cos(\x)}, {2.5 * sin(\x)});
\draw [dotted, line width = 0.5mm,red,domain=275:280] plot ({4.5 * cos(\x)}, {2.5 * sin(\x)});
\draw [line width = 0.3mm,black,domain=280:330] plot ({4.5 * cos(\x)}, {2.5 * sin(\x)});

\draw[->,line width=0.3mm,red] (-4.432,0.434) .. controls (-3.8,1.2) and (-5.5,1.4) .. (-6,1.5);
\draw[red] (-6.2,1.5) node {?};

\draw[->,line width = 0.3mm,blue] (3.6862,1.4339) -- (3.3,0.9);
\draw[->,line width = 0.3mm,blue] (3.4472,1.6070) -- (3.02,1.2);
\draw[->,line width = 0.3mm,blue] (3.182,1.7678) -- (2.8,1.35);
\draw[->,line width = 0.3mm,blue] (2.7705,1.97) -- (2.5,1.55);
\draw[->,line width = 0.3mm,blue] (-3.6862,-1.4339) -- (-3.3,-0.9);
\draw[->,line width = 0.3mm,blue] (-3.4472,-1.6070) -- (-3.02,-1.2);
\draw[->,line width = 0.3mm,blue] (-3.182,-1.7678) -- (-2.8,-1.35);
\draw[->,line width = 0.3mm,blue] (-2.7705,-1.97) -- (-2.5,-1.55);
\draw[->,line width = 0.3mm,blue] (1.4651,2.3638) -- (1.42,2.05);
\draw[->,line width = 0.3mm,blue] (0.7814,2.4620) -- (0.75,2);
\draw[->,line width = 0.3mm,blue] (0.3139,2.4939) -- (0.28,2.1);
\draw[->,line width = 0.3mm,blue] (-0.3139,2.4939) -- (-0.27,2);      
\draw[->,line width = 0.3mm,blue] (-1.4651,-2.3638) -- (-1.42,-2.05);
\draw[->,line width = 0.3mm,blue] (-0.7814,-2.4620) -- (-0.75,-2);
\draw[->,line width = 0.3mm,blue] (-0.3139,-2.4939) -- (-0.28,-2.1);
\draw[->,line width = 0.3mm,blue] (0.3139,-2.4939) -- (0.27,-2);
\draw[->,line width = 0.3mm,blue] (-1.3157,2.33908) -- (-1.27,2.05);
\draw[->,line width = 0.3mm,blue] (-1.9018,2.2658) -- (-1.85,2);
\draw[->,line width = 0.3mm,blue] (-2.7705,1.97) -- (-2.3,1.8);
\draw[->,line width = 0.3mm,blue] (-3.182,1.7678) -- (-2.8,1.6);      
\draw[->,line width = 0.3mm,blue] (-3.5939,1.5045) -- (-3.2,1.4); 
\draw[->,line width = 0.3mm,blue] (-3.8573,1.2876) -- (-3.45,1); 
\draw[->,line width = 0.3mm,blue] (1.3157,-2.33908) -- (1.27,-2.05);
\draw[->,line width = 0.3mm,blue] (1.9018,-2.2658) -- (1.85,-2);
\draw[->,line width = 0.3mm,blue] (2.7705,-1.97) -- (2.3,-1.8);
\draw[->,line width = 0.3mm,blue] (3.182,-1.7678) -- (2.8,-1.6);      
\draw[->,line width = 0.3mm,blue] (3.5939,-1.5045) -- (3.2,-1.4); 
\draw[->,line width = 0.3mm,blue] (3.8573,-1.2876) -- (3.45,-1); 
\end{tikzpicture}
\end{center}
\caption{Representation of the ``partial invariance property'' obtained by application of the results in Theorem~\ref{mainTheorem0}. A priori, we know nothing about the images of the (dotted) red points under \eqref{eq:switchedSystem}. Our goal is to convert this partial invariance property into a global stability property.}
\end{figure}
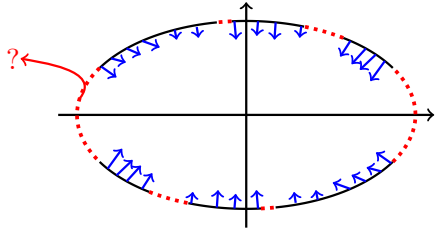

\noindent
Thus, we are left with the following question:

\begin{problem}\label{prob:question}
What can we conclude on the JSR if the invariance property is satisfied by all points, except a set of measure $\tilde{\varepsilon}$?
\end{problem}

\noindent
In the course of Theorem~\ref{thm:mainTheorem01}, we will be able to derive an upper bound by solving the geometric problem of computing the largest scaling of $E_P$ included in the convex hull of the subset of points of $E_P$ that satisfy the invariance property \eqref{eqn:P0}. Indeed, this smaller ellipsoid will satisfy a (relaxed) invariance property \emph{for all its points}, thanks to the following key property of switched linear systems.

\begin{property}\label{property:convpres}
The dynamics given in \eqref{eq:switchedSystem} is convexity-preserving, meaning that for any set of points $X \subset \mathbb{R}^n$, 
\begin{equation*}
\begin{aligned}
&\{Ax: A \in \mathcal{M}, x \in \conv{X}\} \subset \conv{ \{Ax: A\in \mathcal{M}, x\in X \} }.
\end{aligned}
\end{equation*}
\end{property}

\noindent
Of course, for a fixed measure $\tilde{\varepsilon}$, this largest ellipsoid will depend on the distribution of points of $E_P$ that violate the constraint. In order to obtain a guarantee on our upper bound, we will look for the smallest such ellipsoid obtained over all possible sets $V(\omega_N)$ of measure $\tilde{\varepsilon}$. 

\noindent
We start by solving this problem in the particular case where $E_P = \sphere$. In this case, we benefit from the following tool, allowing to explicitly analyse the worst-case distribution.

\bde
We define the \emph{spherical cap} on $\sphere$ for a given hyperplane $c^Tx = k$, as $\mathcal{C}_{c,k} := \{x \in \sphere : c^Tx >k\}$.
\ede

\noindent
We now define the following function which quantifies the largest-inscribed-sphere problem, for a given subset $X \subset \sphere$:
\begin{equation}\label{shrinkage}
\Delta: \left\{
    \begin{split}
    &\wp(\sphere) \to [0,1]\\ 
    &X \mapsto \sup \{r: r\ball \subset \conv(\sphere \setminus X)\}.
    \end{split}
  \right.
\end{equation}
The following proposition tells us that, when the measure of the set $X$ is fixed, $\Delta$ is minized when $X$ is a spherical cap, i.e., the minimal radius $\delta$ of the largest sphere $\delta \sphere$ included in $\sphere \setminus X$ will be reached when $X$ is a spherical cap.

\begin{prop}\label{prop:mainSphericalCap}
Let $\tilde{\varepsilon} \in [0,1]$ and $\mathcal{X}_{\tilde{\varepsilon}} = \{X \subset \sphere: \sigma^{n-1}(X) \leq \tilde{\varepsilon}\}$. Then, the function $\Delta(X)$ attains its minimum over $\mathcal{X}_{\tilde{\varepsilon}}$ for some $X$ which is a spherical cap. We denote by $\delta(\tilde{\varepsilon})$ this minimal value, which takes the following expression: $$\delta(\tilde{\varepsilon}) = \sqrt{1 - I^{-1}(2 \tilde{\varepsilon};\frac{n-1}{2},\frac{1}{2})}.$$
\end{prop}

\noindent
A proof of Proposition~\ref{prop:mainSphericalCap} is given in the Appendix. By homogeneity of Program (10), we have $x \in \tilde{\sphere} \iff -x \in \tilde{\sphere}$, which implies that the minimal $\delta$ will in fact occur when the set of violating points is the union of two symmetric spherical caps, each of measure $\frac{\tilde{\varepsilon}}{2}$.

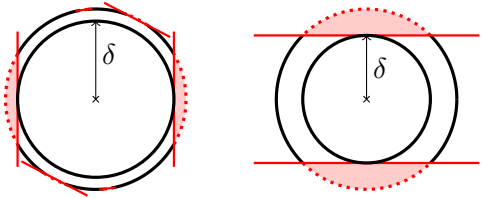
\begin{figure}[H]
\begin{center}
\begin{tikzpicture}[scale=0.6]
\draw (0.07,-0.07) -- (-0.07,0.07);
\draw (0.07,0.07) -- (-0.07,-0.07);
\fill[fill=red!20] (7.414,1.414) arc [start angle=45, end angle=135, radius=2] -- cycle;
\fill[fill=red!20] (4.586,-1.414) arc [start angle=225, end angle=315, radius=2] -- cycle;

\fill[fill=red!20] (1.732,-1) arc [start angle=-30, end angle=30, radius=2] -- cycle;
\fill[fill=red!20] (-1.732,1) arc [start angle=150, end angle=210, radius=2] -- cycle;
\fill[fill=red!20] (1.1472,1.6384) arc [start angle=55, end angle=70, radius=2] -- cycle;
\fill[fill=red!20] (-1.1472,-1.6384) arc [start angle=235, end angle=250, radius=2] -- cycle;
\fill[fill=red!20] (-0.1744,1.9924) arc [start angle=95, end angle=100, radius=2] -- cycle;
\fill[fill=red!20] (0.1744,-1.9924) arc [start angle=275, end angle=280, radius=2] -- cycle;

\draw[dotted, very thick,red,domain=-30:30] plot ({2*cos(\x)}, {2*sin(\x)});
\draw[very thick,black,domain=30:55] plot ({2 * cos(\x)}, {2 * sin(\x)});
\draw[dotted, very thick,red,domain=55:70] plot ({2 * cos(\x)}, {2 * sin(\x)});
\draw[very thick,black,domain=70:95] plot ({2 * cos(\x)}, {2 * sin(\x)});
\draw[dotted, very thick,red,domain=95:100] plot ({2 * cos(\x)}, {2 * sin(\x)});
\draw[very thick,black,domain=100:150] plot ({2 * cos(\x)}, {2 * sin(\x)});
\draw[dotted, very thick,red,domain=150:210] plot ({2 * cos(\x)}, {2 * sin(\x)});
\draw[very thick,black,domain=210:235] plot ({2 * cos(\x)}, {2 * sin(\x)});
\draw[dotted, very thick,red,domain=235:250] plot ({2 * cos(\x)}, {2 * sin(\x)});
\draw[very thick,black,domain=250:275] plot ({2 * cos(\x)}, {2 * sin(\x)});
\draw[dotted, very thick,red,domain=275:280] plot ({2 * cos(\x)}, {2 * sin(\x)});
\draw[very thick,black,domain=280:330] plot ({2 * cos(\x)}, {2 * sin(\x)});

\draw (5.93,-0.07) -- (6.07,0.07);
\draw (5.93,0.07) -- (6.07,-0.07);
\draw [black,very thick,domain=0:45] plot ({6+2*cos(\x)}, {2*sin(\x)});
\draw [dotted, red,very thick,domain=45:135] plot ({6+2*cos(\x)}, {2*sin(\x)});
\draw [black,very thick,domain=135:225] plot ({6+2*cos(\x)}, {2*sin(\x)});
\draw [dotted, red,very thick,domain=225:315] plot ({6+2*cos(\x)}, {2*sin(\x)});
\draw [black,very thick,domain=315:360] plot ({6+2*cos(\x)}, {2*sin(\x)});
\draw[line width = 0.3mm,red] (3.5,1.414) -- (8.5,1.414);
\draw[line width = 0.3mm,red] (3.5,-1.414) -- (8.5,-1.414);

\draw [red, line width = 0.3mm] (1.7321, -1.5) -- (1.7321, 1.5);
\draw [red, line width = 0.3mm] (-1.7321, -1.5) -- (-1.7321, 1.5);
\draw [red, line width = 0.3mm] (1.6472, 1.3781) -- (0.184, 2.1397);
\draw [red, line width = 0.3mm] (-1.6472, -1.3781) -- (-0.184, -2.1397);
\draw [red, line width = 0.3mm] (-0.0743,2.0056) -- (-0.4473,1.9564);
\draw [red, line width = 0.3mm] (0.0743,-2.0056) -- (0.4473,-1.9564);

\draw[->] (6,0) -- (6,1.414);
\draw (6.3,0.7) node {$\delta$};
\draw [black,very thick] (6,0) circle [radius = 1.414cm];
\draw [black,very thick] (0,0) circle [radius = 1.732cm];
\draw[->] (0,0) -- (0,1.732);
\draw (0.3,1) node {$\delta$};
\end{tikzpicture}
\end{center}
\caption{On the left, a general case of the situation where the ellipse in Fig. 1. is a sphere. On the right, case giving minimal $\delta$. The set of points violating the invariance constraint (in red) is the union of two spherical caps.}
\end{figure}
\begin{rem}
When $\varepsilon \geq \frac{1}{m^l}$, we have $\tilde{\varepsilon} \geq 1$ and $\delta(\tilde{\varepsilon}) = 0$: the only upper bound we can give for the JSR is then $+ \infty$.
\end{rem}
\end{subsection}

\begin{subsection}{A global upper bound}
We now introduce our main theorem, Theorem~\ref{thm:mainTheorem01}, which provides a solution to Problem~\ref{prob:question}. In order to use our solution of previous section, developed for the case $E_P = \mathbb{S}$, we will have to relate $E_{P(\omega_N)}$ to $\sphere$. We apply thus a change of coordinates bringing $E_P$ to $\sphere$. Since $P \in \mathcal{S}_{++}^n$, it can be written in its Cholesky form 
\begin{equation}\label{cholesky}
P = L^TL,
\end{equation} 
where $L$ is an upper triangular matrix. Remark that $L$ maps the elements of $E_P$ to $\sphere$. Since the JSR is not changed by similarity transformations, we can pursue our calculations with the matrices obtained after the change of coordinates.
\begin{figure}[H]
\begin{center}
\begin{tikzpicture}[scale=0.4]
\draw (-4-0.1,-0.1) -- (-4+0.1,0.1);
\draw (-4-0.1,0.1) -- (-4+0.1,-0.1);
\draw [dotted, line width = 0.7mm,red,domain=-30:30] plot ({-4 + 4.5 * cos(\x)}, {2.5 * sin(\x)});
\draw [line width = 0.5mm,black,domain=30:55] plot ({-4 + 4.5 * cos(\x)}, {2.5 * sin(\x)});
\draw [dotted, line width = 0.7mm,red,domain=55:70] plot ({-4 + 4.5 * cos(\x)}, {2.5 * sin(\x)});
\draw [line width = 0.5mm,black,domain=70:95] plot ({-4 + 4.5 * cos(\x)}, {2.5 * sin(\x)});
\draw [dotted, line width = 0.7mm,red,domain=95:100] plot ({-4 + 4.5 * cos(\x)}, {2.5 * sin(\x)});
\draw [line width = 0.5mm,black,domain=100:150] plot ({-4 + 4.5 * cos(\x)}, {2.5 * sin(\x)});
\draw [dotted, line width = 0.7mm,red,domain=150:210] plot ({-4 + 4.5 * cos(\x)}, {2.5 * sin(\x)});
\draw [line width = 0.5mm,black,domain=210:235] plot ({-4 + 4.5 * cos(\x)}, {2.5 * sin(\x)});
\draw [dotted, line width = 0.7mm,red,domain=235:250] plot ({-4 + 4.5 * cos(\x)}, {2.5 * sin(\x)});
\draw [line width = 0.5mm,black,domain=250:275] plot ({-4 + 4.5 * cos(\x)}, {2.5 * sin(\x)});
\draw [dotted, line width = 0.7mm,red,domain=275:280] plot ({-4 + 4.5 * cos(\x)}, {2.5 * sin(\x)});
\draw [line width = 0.5mm,black,domain=280:330] plot ({-4 + 4.5 * cos(\x)}, {2.5 * sin(\x)});

\draw (5.9,-0.1) -- (6.1,0.1);
\draw (5.9,0.1) -- (6.1,-0.1);
\draw [dotted, line width = 0.7mm,red,domain=-30:30] plot ({6 + 2 * cos(\x)}, {2 * sin(\x)});
\draw [line width = 0.5mm,black,domain=30:55] plot ({6 + 2 * cos(\x)}, {2 * sin(\x)});
\draw [dotted, line width = 0.7mm,red,domain=55:70] plot ({6 + 2 * cos(\x)}, {2 * sin(\x)});
\draw [line width = 0.5mm,black,domain=70:95] plot ({6 + 2 * cos(\x)}, {2 * sin(\x)});
\draw [dotted, line width = 0.7mm,red,domain=95:100] plot ({6 + 2 * cos(\x)}, {2 * sin(\x)});
\draw [line width = 0.5mm,black,domain=100:150] plot ({6 + 2 * cos(\x)}, {2 * sin(\x)});
\draw [dotted, line width = 0.7mm,red,domain=150:210] plot ({6 + 2 * cos(\x)}, {2 * sin(\x)});
\draw [line width = 0.5mm,black,domain=210:235] plot ({6 + 2 * cos(\x)}, {2 * sin(\x)});
\draw [dotted, line width = 0.7mm,red,domain=235:250] plot ({6 + 2 * cos(\x)}, {2 * sin(\x)});
\draw [line width = 0.5mm,black,domain=250:275] plot ({6 + 2 * cos(\x)}, {2 * sin(\x)});
\draw [dotted, line width = 0.7mm,red,domain=275:280] plot ({6 + 2 * cos(\x)}, {2 * sin(\x)});
\draw [line width = 0.5mm,black,domain=280:330] plot ({6 + 2 * cos(\x)}, {2 * sin(\x)});
\draw[->] (0,2) .. controls (1,3) and (3,3) .. (4,2);
\draw[->] (4,-2) .. controls (3,-3) and (1,-3) .. (0,-2); 
\draw (2,3.3) node {$L$};
\draw (2,-3.3) node {$L^{-1}$};
\end{tikzpicture}
\end{center}
\caption{Change of coordinates to bring our problem back to the case of the unit sphere.}
\end{figure}
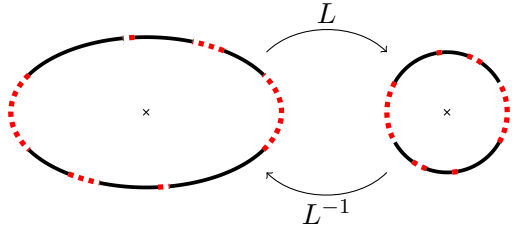
\btheos\label{thm:mainTheorem01}
Let $\gamma^* \in \mathbb{R}_{> 0}$. Consider a set of matrices $\mathcal{M}$, and a matrix $P \succ 0$ optimal solution of $\Opt(\omega_N)$, satisfying Equation \eqref{eqn:P0} for some $\tilde{\sphere} \subset \sphere$ where $\sigma^{n-1}(\tilde{\sphere}) \leq \tilde{\varepsilon}$. Then, we have 
\begin{equation}\label{eq:ub1}
\rho(\mathcal{M}) \leq \frac{\gamma^{*}}{\sqrt[l]{\delta \left( \frac{\tilde{\varepsilon} \kappa(P)}{2} \right)}}
\end{equation} 
with $\kappa(P) = \sqrt{\frac{\det(P)}{\lambda_{\min}(P)^n}}$ and where $\delta(\cdot{})$ is given by Proposition~\ref{prop:mainSphericalCap}.
\etheos
\HACKPROOF
\ \newline
\textit{i)} Since we have seen in the previous section a technique to solve the spherical case, we first bring our problem to the spherical case. To do so, we perform the change of coordinates defined as in \eqref{cholesky} by $L \in \mathcal{L}(\mathbb{R}^n)$ which maps the ellipsoid $E_P$ to the sphere $\sphere$. By defining \mbox{$\bar{A}_{j_i}=  L A_{j_i} L^{-1}$,} and $\mathbf{\bar{A}_j} = \bar{A}_{j_l} \bar{A}_{j_{l-1}} \dots \bar{A}_{j_1} $, Equation \eqref{eqn:P0} becomes:
\begin{equation}\label{eq:assumption3}
\hspace{-0.3cm}(\mathbf{\bar{A}_j} x)^T \mathbf{\bar{A}_j} x \leq ({\gamma^*})^{2l} x^T x, \, \forall x \in L(\sphere \setminus \tilde{\sphere}), \forall \mathbf{j} \in M^l.
\end{equation}
By using the homogeneity of the dynamics, we have:
\begin{equation*}
\begin{aligned}
&(\mathbf{\bar{A}_j} x)^T \mathbf{\bar{A}_j} x \leq ({\gamma^*})^{2l} x^T x,\, \forall x \in L(\sphere \setminus \tilde{\sphere}) \implies (\mathbf{\bar{A}_j} x)^T \mathbf{\bar{A}_j} x \leq ({\gamma^*})^{2l} x^T x,\, \forall x \in \Pi_\sphere \left( L(\sphere \setminus \tilde{\sphere}) \right),
\end{aligned}
\end{equation*}
and therefore we can rewrite \eqref{eq:assumption3} as:
\begin{equation}\label{eq:assumption4}
(\mathbf{\bar{A}_j} x)^T \mathbf{\bar{A}_j} x \leq ({\gamma^*})^{2l} x^T x, \forall x \in \sphere \setminus \Pi_\sphere(L(\tilde{\sphere}), \forall \mathbf{j} \in M^l.
\end{equation}

\noindent
\textit{ii)} We now show how to relate $\sigma^{n-1}(\proj_\sphere(L(\tilde{\sphere})))$ to $\sigma^{n-1}(\tilde{\sphere})$, the measure of the violating set in the initial coordinates. Consider $\sphere^{\tilde{\sphere}}$, the sector of $\ball$ defined by $\tilde{\sphere}$. We denote $C:= L(\tilde{\sphere})$ and $C':=\proj_\sphere(L(\tilde{\sphere}))$. Thus, we have $\sphere^{C'} \subset \frac{1}{ \lambda_{\min (L^)}} \sphere^C$. This leads to\footnote{Recall that $\lambda(S)$ is the Lebesgue measure of $S$, and the spherical measure of any set $C \subset \mathbb{S}$ is given by $\sigma^{n-1}(C)=\lambda \left( \mathbb{S}^C \right)$.}:
$$\sigma^{n-1}(C') = \lambda \left( \sphere^{C'} \right) \leq \lambda \left( \frac{1}{ \lambda_{\min (L)}} L(\sphere^{\tilde{\sphere}}) \right).$$ Then, the following holds: 
\begin{eqnarray}
\nonumber\sigma^{n-1}(C') &\leq& \frac{\lambda \left( L(\sphere^{\tilde{\sphere}}) \right)}{\lambda_{\min}(L)^n} \\
\label{eqn:lt} &=&\frac{|\det(L)|}{\lambda_{\min}(L)^n}\lambda \left( \sphere^{\tilde{\sphere}} \right) \\
\label{eqn:map1} &=&\sqrt{\frac{\det(P)}{\lambda_{\min}(P)^n}}\sigma^{n-1}(\tilde{\sphere})
\end{eqnarray}
where \eqref{eqn:lt} follows from the fact that $\lambda(Q(X)) = |\det(Q)| \lambda(X)$, for any set $X \subset \mathbb{R}^n$ and $Q \in \mathcal{L}(\mathbb{R}^n)$ (see e.g. \cite{rudin}). Hence, we have
\begin{equation}\label{eqn:contraction}
(\mathbf{\bar{A}_j} x)^T \mathbf{\bar{A}_j} x \leq (\gamma^{*})^{2l} x^Tx, \forall x \in \sphere \setminus \sphere', \forall \mathbf{j} \in M^l,
\end{equation}
with $\sphere' = \Pi_\sphere(L(\tilde{\sphere})$ and $\sigma^{n-1}(\sphere') \leq \sqrt{\frac{\det(P)}{\lambda_{\min}(P)^n}}\sigma^{n-1}(\tilde{\sphere})= \kappa(P) \tilde{\varepsilon}$.

\noindent
\textit{iii)} For any such given set $\sphere'$, we look for the largest sphere included in $\convhull (\sphere \setminus \sphere')$. By homogeneity of the system, this sphere is centered at the origin, and we denote by $\alpha$ its radius. By \eqref{eqn:contraction}, $l$-traces initialized on $\sphere \setminus \sphere'$ will be in $(\gamma^*)^l \ball$: $$\mathbf{\bar{A}_j} \left( \sphere \setminus \sphere' \right) \subset (\gamma^*)^l \ball, \, \forall \mathbf{j} \in M^l.$$
Now, combining with Property~\ref{property:convpres}, we have: $$\mathbf{\bar{A}_j} \left( \convhull(\sphere \setminus \sphere') \right) \subset \convhull( \mathbf{\bar{A}_j} (\sphere \setminus \sphere')) \subset (\gamma^*)^l \ball, \forall \mathbf{j} \in M^l.$$ Since $\alpha \sphere \subset \convhull (\sphere \setminus \sphere')$, then $\forall \mathbf{j} \in M^l$, $\mathbf{\bar{A}_j} \left( \alpha \sphere \right) = \alpha \mathbf{\bar{A}_j} \left( \sphere \right) \subset \convhull \left( \mathbf{\bar{A}_j} (\sphere \setminus \sphere') \right) \subset (\gamma^*)^l \ball$, which implies that 
\begin{equation}\label{eq:part3}
\mathbf{\bar{A}_j} (\sphere) \subset \frac{(\gamma^*)^l}{\alpha} \ball.
\end{equation}

\noindent
\textit{iv)} Summarizing, since we know that the set $\tilde{\mathbb{S}}$ is symmetric w.r.t. the origin, by Proposition 1, we have that $\alpha \geq \delta(\frac{\tilde{\varepsilon} \kappa(P)}{2})$. Finally, by homogeneity of our system and the fact that the JSR is invariant under similarity transformations, Equation \eqref{eq:part3} implies $\rho(\mathcal{M}^l) \leq \frac{(\gamma^*)^l}{\delta(\frac{\tilde{\varepsilon} \kappa(P)}{2})}$, hence $\rho(\mathcal{M}) \leq \frac{\gamma^{*}}{\sqrt[l]{\delta \left( \frac{\tilde{\varepsilon} \kappa(P)}{2} \right)}}$.
\HACKENDPROOF

\begin{rem}\label{otherUb}
There is no conservatism in multiplying $\varepsilon$ by $m^l$, as in the worst case this can happen: if $\varepsilon = 1/m^l$, Theorem 4 does not rule out the pathological case where not a single point satisfies Equation (7) for all $\mathbf{A} \in \mathcal{M}^l$, and thus $\delta$ must be equal to $0$ since all points might be violating the constraint. However, the multiplication by $\kappa(P)$ is conservative if $P$ has different eigenvalues (this bound is then exactly reached only at a single point on the ellipsoid). We can then, instead of deriving an upper bound on the size of the set of points that violate the constraint, look at a lower bound on the size of the set of points that satisfy the constraint. Taking the complement of this latter set gives another upper bound on the size of the set of violating points. By a similar reasoning as the one conducted above, this second upper bound will be equal to $1 - (1-\varepsilon m^l) \sqrt{\frac{\det(P)}{\lambda_{\max}(P)^n}}$. \\
This provides an alternative upper bound, which can be used if the initial upper bound \eqref{eq:ub1} is infinite, or weaker.
\end{rem}
\end{subsection}

\begin{subsection}{Main Theorem}
We are now ready to prove our main theorem by putting together all the above pieces.

\btheos\label{thm:mainTheorem}
Consider an $n$-dimensional switched linear system as in \eqref{eq:switchedSystem} and a uniform random sampling $\omega_N \subset Z_l$, where $N \geq \frac{n(n+1)}{2}+1$. For any $\eta > 0$, let $\gamma^{*}(\omega_N) $ be the optimal solution to \eqref{eq:lowerbound}. Then, for any given $\beta \in [0,1)$, with probability at least $\beta$ we have:
$$\rho \leq \frac{\gamma^{*}(\omega_N) (1+ \eta)}{\sqrt[l]{\delta(\beta, \omega_N)}},$$
where $\delta(\beta, \omega_N) = \sqrt{1 - I^{-1}(\varepsilon(\beta,N) m^l \kappa(P); \frac{n-1}{2}, \frac{1}{2})}$ satisfies $\lim_{N \to \infty}\delta(\beta, \omega_N) = 1$ with probability $1$.
\etheos
\HACKPROOF
Let us consider $\gamma^*(\omega_N)$ and $P$ as in Equation \eqref{eqn:campiOpt03}. Then, by taking $\varepsilon:=\varepsilon(\beta,N)$ such that $\beta(\varepsilon,N)=\beta$ in Corollary~\ref{cor:gettingRidOfm}, we have 
\begin{equation*} 
(\mathbf{A_j} x)^T P \mathbf{A_j} x \leq  \left( (\gamma^{*}(1+\eta) \right)^{2l} x^T P x, \forall x \in \sphere \setminus \tilde{\sphere}, \forall \mathbf{j} \in M^l
\end{equation*} 
with $\tilde{\sphere}$ the projection of $V$ on $\sphere$, and $\sigma^{n-1}(\tilde{\sphere}) \leq \varepsilon m^l$. Then by Theorem \ref{thm:mainTheorem01}, we can compute $\delta(\beta, \omega_N) =\delta(\varepsilon'(\beta,N))$, where
\begin{equation}\label{eqn:eps2}
\varepsilon'(\beta, N) = \frac{1}{2} \varepsilon(\beta,N) m^l \kappa(P) 
\end{equation} 
such that with probability at least $\beta$ we have: $$\rho \leq \frac{\gamma^{*}(\omega_N) (1 + \eta)}{\sqrt[l]{\delta(\beta, \omega_N)}},$$ which completes the proof of the first part of the theorem. 

\vspace{0.3cm}

\noindent
Let us prove now that $\lim_{N \to \infty} \delta(\beta, \omega_N) = 1$ with probability $1$. We recall that $$\delta(\beta, \omega_N) = \delta \left( \varepsilon(\beta, \omega_N) m^l \kappa(P(\omega_N)) \right).$$ 

\noindent
We start by showing that $\kappa \left( P(\omega_N) \right)$  is uniformly bounded in $N$. The optimization problem $\Opt(\omega_N)$ given in \eqref{eqn:campiOpt03}, with $(\omega_N)$ replaced by $(Z_l)$ and $(1 + \eta)$ replaced by $(1+\frac{\eta}{2})$ is strictly feasible for any positive parameter $\eta$. It then admits a finite optimal value $K$ for some solution $P_{\eta/2}$. Note that, $\lim_{N \to \infty} \gamma^{*}(\omega_N)= \gamma^{*}(Z_l)$ with probability $1$. Thus, for large enough $N$, \mbox{$\gamma^{*}(\omega_N)(1+\eta) > \gamma^{*}(Z_l)(1+\frac{\eta}{2})$.} This also means that, for large enough $N$, $\Opt(\omega_N)$ admits $P_{\eta/2}$ as a feasible solution and thus the optimal value of $\Opt(\omega_N)$ is upper-bounded by $K.$ In other words, for $N$ large enough, \mbox{$\lambda_{\max}(P({\omega_N})) \leq K$.} Moreover, since  $\lambda_{\min}(P(\omega_N))\geq 1$ (by $P \succeq I $), we also have \mbox{$\det(P(\omega_N)) \geq 1$,} which means that
\begin{equation}\label{kappa}
\kappa \left( P(\omega_N) \right) = \sqrt{\frac{\det(P(\omega_N))}{\lambda_{\min}(P(\omega_N))^n}} \leq \sqrt{K^n}.
\end{equation}
We next show that $\lim_{N \to \infty} \varepsilon(\beta, N) = 0$ for any fixed $\beta \in [0,1)$. Recall that $\varepsilon(\beta, N)$ is implicitely defined by
\begin{eqnarray}
\nonumber 1-\beta &=& \sum_{j=0}^d {{N}\choose{j}} \varepsilon^j (1-\varepsilon)^{N-j}\\
                 &\leq& (d+1)N^d (1-\varepsilon)^{N-d}.\label{eqn:beta}
\end{eqnarray}

\noindent 
We prove $\lim_{N \to \infty} \varepsilon(\beta, N) = 0$ by contradiction. Assume that $\lim_{N \to \infty} \varepsilon(\beta, N) \not= 0$. This means that, there exists some $c > 0$ such that $\varepsilon(\beta, N) > c$ infinitely often. Then, consider the subsequence $N_k$ such that $\forall k$, $\varepsilon(\beta, N_k) > c$. Then, by \eqref{eqn:beta} we have for any $k \in \mathbb{N}$:
\begin{equation*}
1-\beta \leq  (d+1)N_k^d (1-\varepsilon)^{N_k-d}\hspace{-0.4mm} \leq \hspace{-0.7mm}(d+1)N_k^d (1-c)^{N_k-d}. 
\end{equation*}
Note that $\lim_{k \to +\infty}(d+1)N_k^d (1-c)^{N_k-d} = 0$, which implies that there exists a $k'$ such that:
$$(d+1)N_{k'}^d (1-c)^{N_k'-d} < 1 - \beta,$$ which is a contradiction. Therefore, we must have  $\lim_{N \to \infty} \varepsilon (\beta, N) = 0$. Putting this together with \eqref{kappa}, we get: $\lim_{N \to \infty} m^l \kappa(P(\omega_N)) \varepsilon(\beta, \omega_N) = 0$. By the continuity of the function $\delta$ this also implies: $\lim_{N \to \infty} \delta \left( \varepsilon(\beta, \omega_N) m^l \kappa(P(\omega_N)) \right) = 1.$
\HACKENDPROOF
\end{subsection}

\section{Experimental Results}\label{sec:experiments}
\begin{subsection}{Algorithm and experimental results}
Theorem~\ref{thm:lowerbound} and Theorem~\ref{thm:mainTheorem} give us a straightforward algorithm which is summarized at the end of Section 2. In its first part, we look for $\gamma^*$ by bisection on an interval $[0,U]$ (for the value of $U$, take, e.g. the maximum value of $||x_{k+l}||$ among the observations made). For a fixed desired accuracy $\alpha$ on that bisection, we solve a feasibility problem (of polynomial complexity in the number of constraints) at most $\lceil \log_2(U/\alpha) \rceil$ times. In our experiments we took $\alpha=10^{-3}$. Once the result of the bisection is obtained, we solve $\Opt(\omega_N)$. In practice, the parameter $\eta$ in $\Opt(\omega_N)$ can be put to zero, as it is included in $\alpha$. Finally, we get $\delta$ by using the expression given in Theorem~\ref{thm:mainTheorem01}. All these computations are also of polynomial complexity.

\vspace{0.3cm}

\noindent
We illustrate our technique on a $4$-dimensional switched system with $6$ modes. We fix the confidence level, \mbox{$\beta = 0.95$}, and compute the lower and upper bounds on the JSR for $N:=20+200k,\, k \in\{0, \ldots, 29\}$, according to Theorem~\ref{thm:lowerbound} and Theorem~\ref{thm:mainTheorem}, respectively. We take the average performance of our algorithm over $10$ different runs. Fig.~\ref{fig:exp1} shows the evolution, as $N$ increases, of the upper and lower bounds for various values of trace length $l$. To further demonstrate the practical performance of our technique, we also provide the true value of the JSR approximated by the JSR Toolbox \cite{jsrtoolbox} for this system, which turns out to be $0.918 \pm 0.001$. We observed that the performance of the upper bound is much better for traces of length $1$, while for the lower bound, we benefit by considering traces of higher length. While it is expected that longer traces improve the accuracy, the decreasing performance for the upper bound comes from the fact that many more points are needed for larger traces, because the probability space to be sampled is larger. In our example, our first upper bound smaller than $1$ (that is, being a stability guarantee) was obtained for $N=5820$.

\begin{figure}[h!]
\begin{center}
\includegraphics[scale=0.6]{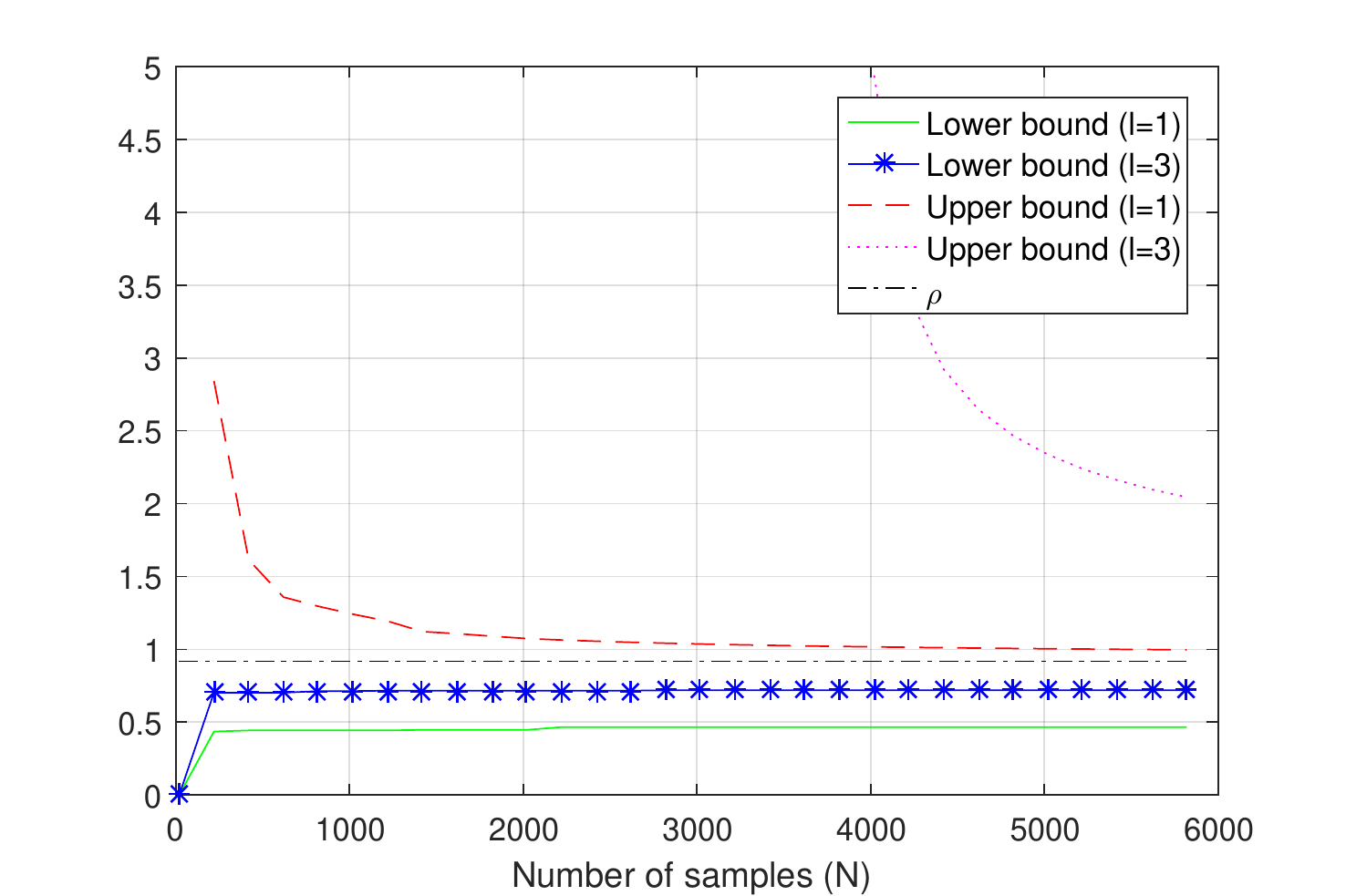}
\end{center}
\caption{Evolution of the upper and lower bounds, for various lengths of traces, with the number of samples.}
\label{fig:exp1}
\end{figure}

\noindent
In Fig.~\ref{fig:exp2}, we compare the upper bound we obtained with the upper bound given by the (white box) JSR Toolbox, for different values of $n$ and $m$.
\begin{figure}[h!]
\begin{center}
\includegraphics[scale=0.6]{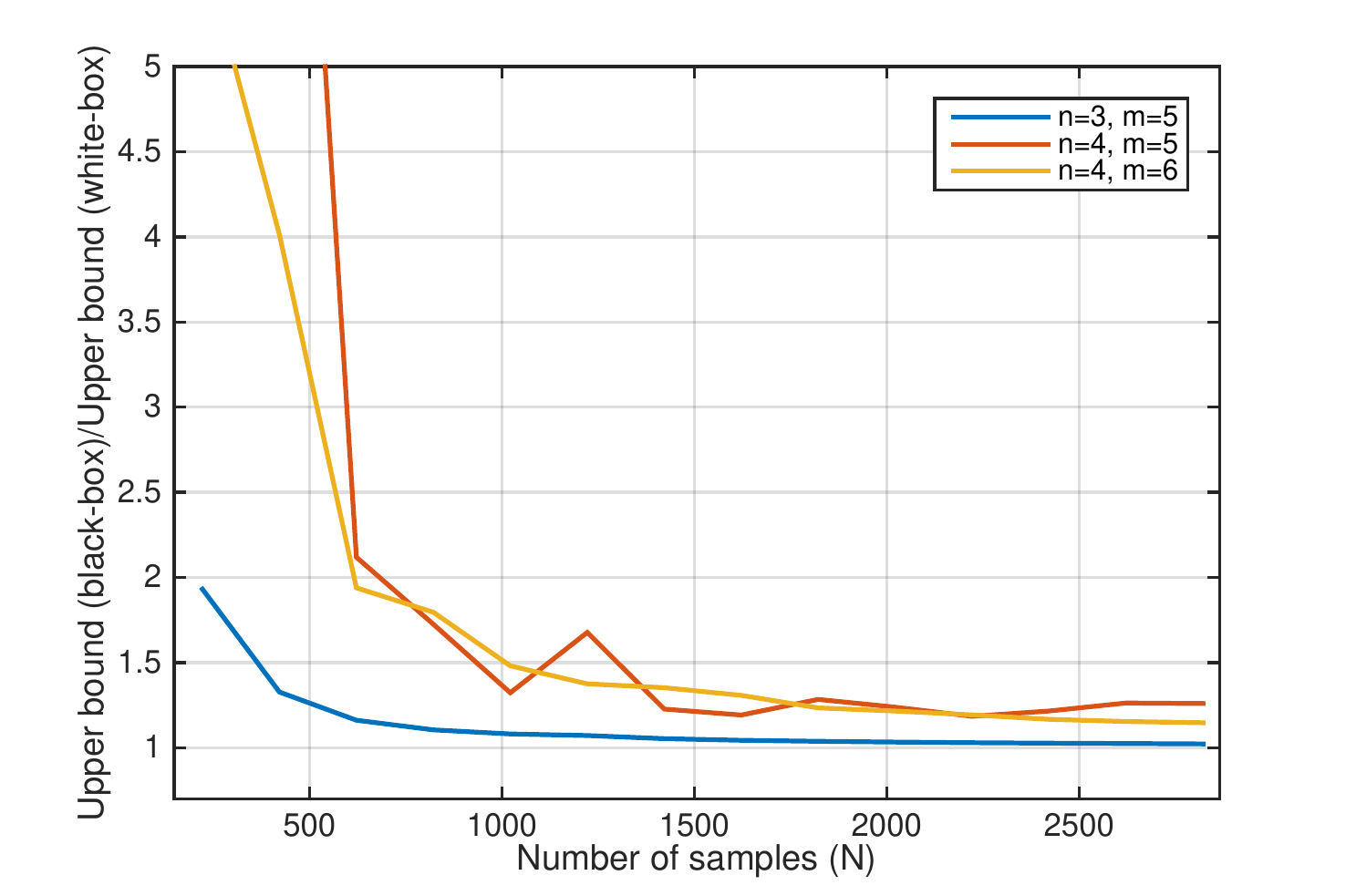}
\end{center}
\caption{Convergence of our upper bound when the number of samples increases, for several values of $n$ and $m$ and $l=1$. The values plotted are the ratios between our upper bound and the true value computed by the JSR Toolbox.}
\label{fig:exp2}
\end{figure}

\noindent
Note that the speed of convergence of all the quantities considered decreases when the dimension of the system increases. We nevertheless observe convergence of the upper bound to $\rho(\mathcal{M})$, and convergence of the lower bound to $\frac{\rho(\mathcal{M})}{\sqrt[l]{n}}$. The gap between these two limits is $\frac{\rho}{\sqrt[l]{n}}$ as predicted by Theorem~\ref{thm:mainTheorem}. This gap could be improved by considering a more general class of common Lyapunov functions, such as those that can be described by sum-of-squares polynomials \cite{sosLyap}. We leave this for future work.

\vspace{0.3cm}

\noindent
To illustrate the accuracy of our confidence level $\beta$, we randomly generate $10,000$ test cases with systems of dimension between $2$ and $7$, number of modes between $2$ and $6$, and size of samples $N$ between $30$ and $1000$. We take $\beta = 0.95$ and we check if the upper bound computed by our technique is greater than the true value computed by the JSR Toolbox for the system. We get $9921$ positive tests, out of $10,000$, which gives us a correctness of $0.9921$ for the upper bound computed. Note that, this probability is significantly above the provided $\beta$. This is expected, since our techniques are based on worst-case analysis and thus fairly conservative.

\vspace{0.3cm}

\noindent
Finally, Fig. \ref{fig:exp3} shows the evolution of the function $\delta$ with the number of samples, for different values of $n$, at $m$ and $l$ fixed.

\begin{figure}[h!]\label{delta}
\begin{center}
\includegraphics[scale=0.5]{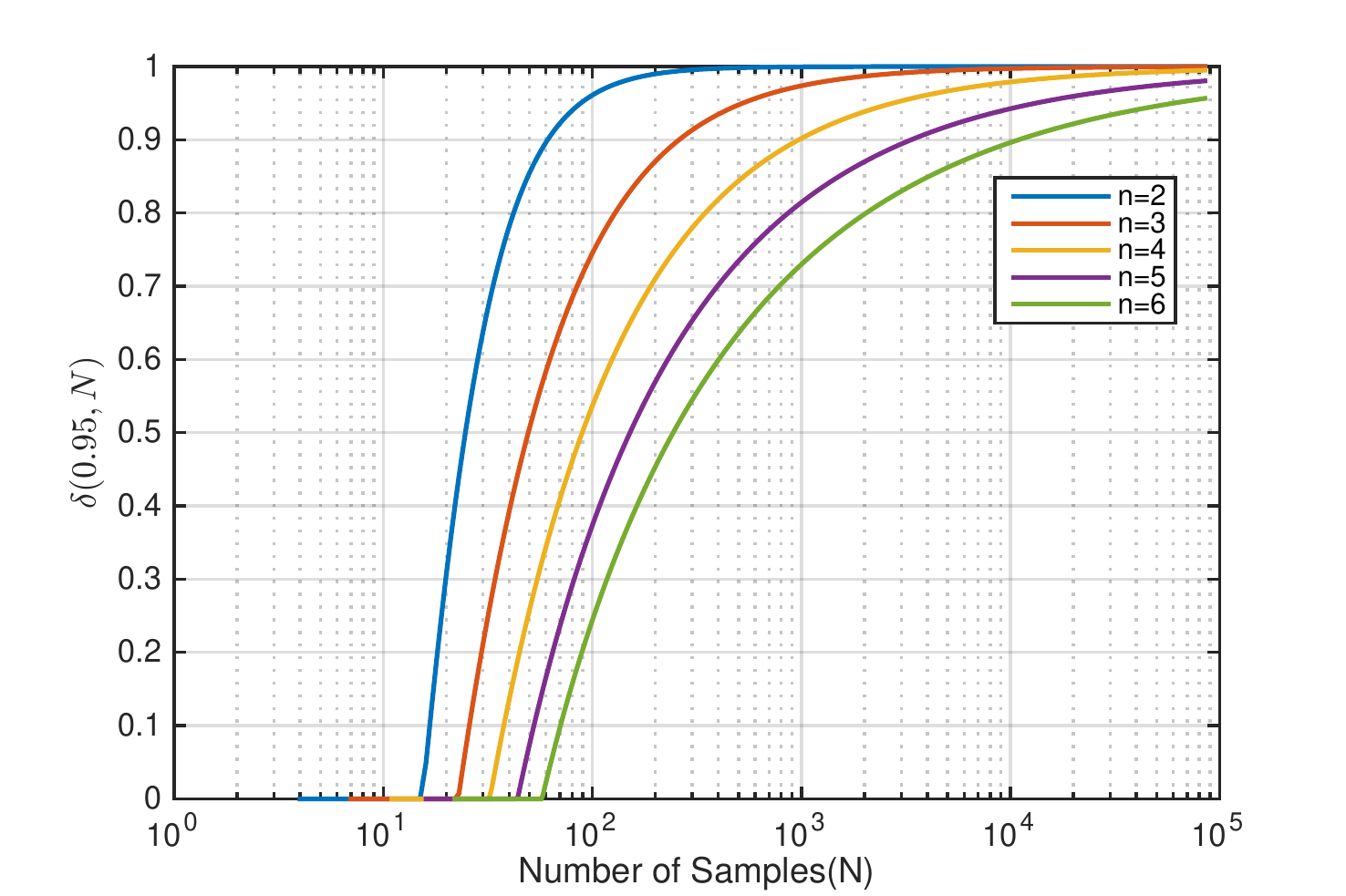}
\end{center}
\caption{Average behavior of $\delta$ as a function of $N$ for different values of $n$, with fixed $m=4$ and $l=1$.}
\label{fig:exp3}
\end{figure}

\end{subsection}

\begin{subsection}{Application to Networked Control System}\label{networkedEx}
We consider here a linear time-invariant control system given as \mbox{$x_{k+1}=Ax_k+Bu_k$}, with control law of the form $u_k = Kx_k$. Matrices $A$, $B$ and $K$ are unknown. The open-loop system is unstable with eigenvalues at $\{0.45, 1.1\}$. The controller stabilizes the system by bringing its eigenvalues to $\{0.8, -0.7\}$.

\vspace{0.3cm}

\noindent
The control input is transmitted over a wireless communication channel that is utilized by $\ell$ users, including the controller. The channel has two modes: a contention access mode, where the users can only send their message if the channel is ``idle'' with carrier-sense multiple access with collision avoidance (CSMA/CA); and a collusion free mode, where each user has guaranteed time slots, during which there is no packet loss.

\vspace{0.3cm}

\noindent
More precisely, the communication between the users and the recipients is performed based on the IEEE 802.15.4 MAC layer protocol \cite{macLayer}, which is used in some of the proposed standards for control over wireless, e.g., WirelessHART \cite{wirelessHart}. This protocol integrates both contention based slots and guaranteed slots. In this example, we consider a beacon-enabled version of the MAC protocol. A centralized control user periodically synchronizes and configures all the users. This control period is named Beacon Interval, and is divided into two subintervals: an active and an inactive period. The active period is itself divided into 6 slots. The first 2 slots constitute the contention access period (CAP), and the next 4 slots constitute the collusion free period (CFP). In our example, the third and fourth slots are designated for the controller, while the fifth and sixth slots are allocated to the other users. Finally, during the inactive period, all users enter a low-power mode to save energy. We illustrate the overall structure of this communication protocol in Fig.~\ref{comExample}. We now want to decide whether the resulting closed-loop networked control system is stable by simulating it starting from different initial conditions. 

\begin{figure}[h!]
\begin{center}
\includegraphics[scale=1.1]{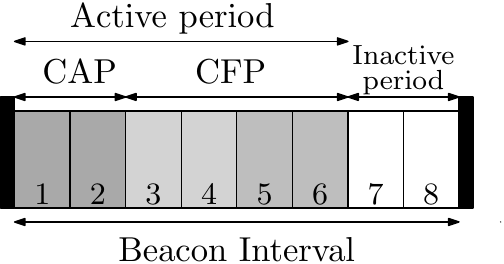}
\end{center}
\caption{The time allocation structure of the modified IEEE 802.15.4 MAC layer.}
\label{comExample}
\end{figure}

\noindent
Note that, the closed-loop dynamics of the underlying system when the controller is active is $A_c = A+BK$. Then, we can model the overall networked control system by the switched linear system $x_{k+1} = \bar{A}x_k$, where $\bar{A} \in \mathcal{M}$ and $\mathcal{M} = \{A^2A_c^2A^4, A_cAA_c^2A^4,AA_c^3A^4, A_c^4A^4\}$. Each element of $\mathcal{M}$ corresponds to a different utilization of the CFP by the users. For example, the mode defined by $A_cAA_c^2A^4$ is active when the first slot in the CFP is assigned to the controller and the second slot is assigned to the other users. We assume that all of the users using the channel have an equal probability of being assigned to a time slot during the CFP. Therefore, the probability of each mode in $\mathcal{M}$ being active is $\left\{\frac{1}{(\ell-1)^2}, \frac{1}{\ell(\ell-1)}, \frac{1}{(\ell-1)\ell}, \frac{1}{\ell^2}\right\}$. Hence, we make use of Remark \ref{rem:probGeneralize} and update our bounds accordingly. Fig.~\ref{fig:networks} shows the computed bounds. As can be seen, approximately after 500 samples, the upper bound on the JSR drops below $1$, which lets us decide that the given closed-loop networked control system is stable, with probability $0.95$.

\begin{figure}[h!]
\begin{center}
\includegraphics[scale=0.6]{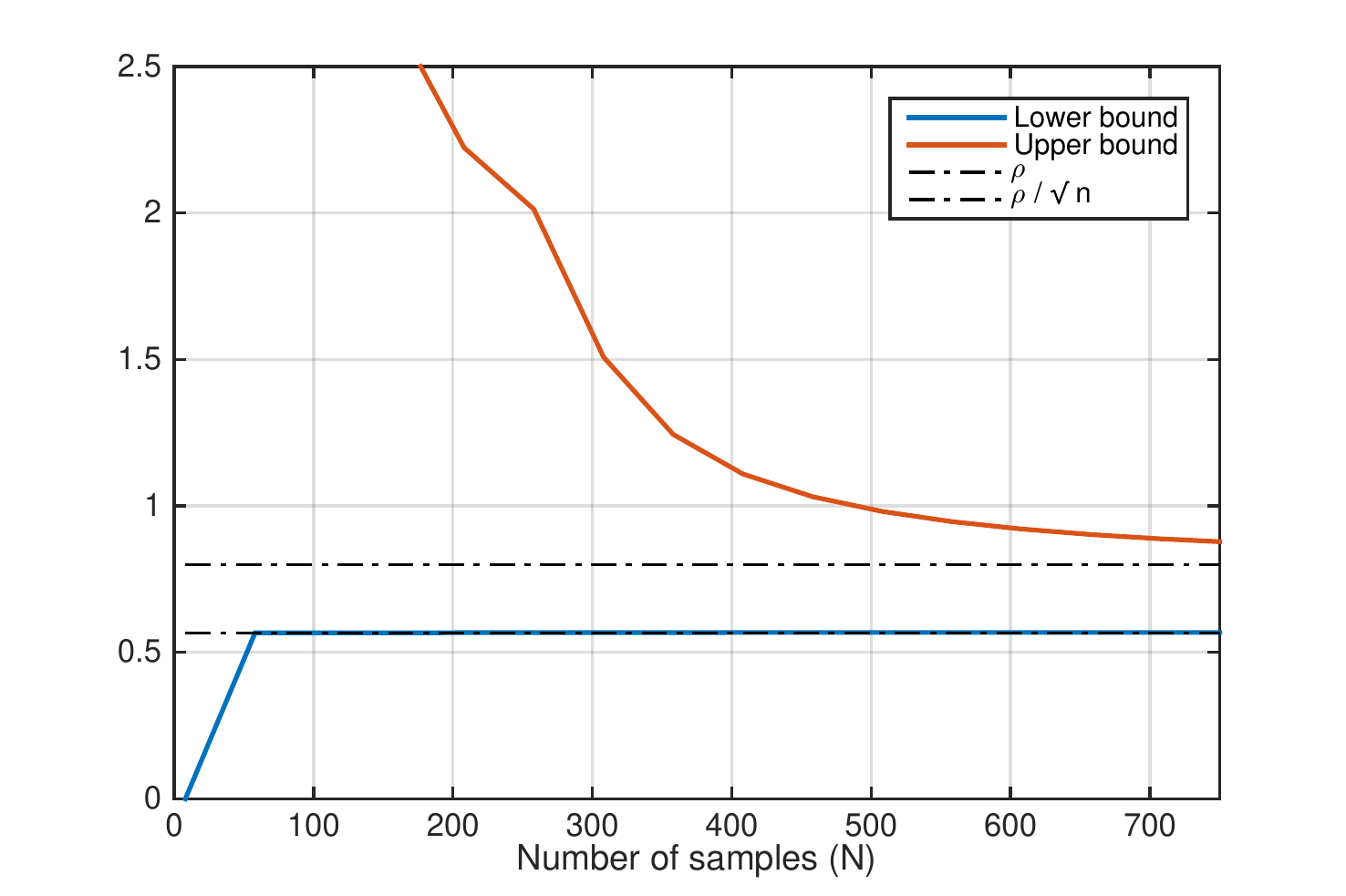}
\end{center}
\caption{The evolution of the computed upper and lower bounds on the JSR with respect to the number of simulations collected from the networked control system.}
\label{fig:networks}
\end{figure}

\end{subsection}

\section{Conclusions}\label{sec:conclusions}
In this paper, we investigated the question of how one can conclude stability of a dynamical system when a model is not available and, instead, we can only observe the evolution of the state of the system. Our goal was to understand how the observation of well-behaved trajectories \emph{intrinsically} implies stability of a system. 

\vspace{0.3cm}

\noindent
As expected, it is not surprising that we need some standing assumptions on the system, in order to allow for any sort of nontrivial stability certificate solely from a finite number of observations. 

\vspace{0.3cm}

\noindent
The novelty of our contribution is twofold: First, we used as standing assumption that the unknown system can be modeled by a switched linear system. This assumption covers a wide range of systems of interest, and to our knowledge no such ``black-box'' result has been available so far on switched systems. Second, we applied powerful techniques from chance-constrained optimization. Their application was far from obvious, though, and relied on geometric properties of linear switched systems.

\vspace{0.3cm}

\noindent
We leveraged the concept of `$l$-step CQF', and showed that it allows to reach arbitrary precision for our black-box technique. In the switched systems literature, there are other well-known techniques for refining this precision for the white-box problem. For instance, one can replace the LMIs in Theorem \ref{thm:cqlf} by Sum-Of-Squares (SOS) constraints; see \cite{parrilo-jadbabaie} or \cite[Theorem 2.16]{jungers_lncis}. Although $l$-step CQFs seem better suited for our purpose, we leave for further work a more systematic analysis of the behaviour of the different refining techniques.

\vspace{0.3cm}

\noindent
Notice that, our algorithm can also be used in the white-box framework and becomes then a randomized algorithm to evaluate the JSR of a known system.

\vspace{0.3cm}

\noindent
In our view, the stability-like guarantees obtained are powerful, in view of the hardness of the general problem. In the future, we plan to investigate how to generalize our results to more complex models of realistic systems.

\section*{Appendix: a Few Results on Spherical Caps}\label{appendix}
Before proceeding to the proof of Proposition~\ref{prop:mainSphericalCap}, we first introduce some necessary definitions and related background on spherical caps. We recall that a \emph{spherical cap} on $\sphere$ for a given hyperplane $c^Tx = k$ is defined by $\mathcal{C}_{c,k} := \{x \in \sphere : c^Tx >k\}$.

\begin{rem}\label{lemma:muMonotone}
Consider the spherical caps $\mathcal{C}_{c, k_1}$ and $\mathcal{C}_{c, k_2}$ such that $k_1 > k_2$, then we have:
$$\sigma^{n-1}(\mathcal{C}_{c,k_1}) < \sigma^{n-1}(\mathcal{C}_{c,k_2}).$$
\end{rem}

\begin{rem}\label{prop:distance}
The distance between the point $x=0$ and the hyperplane $c^Tx = k$ is $\frac{|k|}{\|c\|}$.
\end{rem}

\noindent
We recall the definition given in Section 4 of the function $\Delta$, which quantifies the largest-inscribed-sphere problem, for a given subset $X \subset \sphere.$
\begin{equation*}
\Delta: \left\{
    \begin{split}
    &\wp(\sphere) \to [0,1]\\ 
    &X \mapsto \sup \{r: r\ball \subset \conv(\sphere \setminus X)\}.
    \end{split}
  \right.
\end{equation*}

\begin{lem} \label{lemma:delta2} 
Consider the spherical cap $\mathcal{C}_{c,k}$. We have:
$$\Delta(\mathcal{C}_{c,k}) = \min\left(1, \frac{|k|}{\|c\|}\right).$$
\end{lem}

\HACKPROOF
Note that: $$\conv(\sphere \setminus X)= \left\{x \in \ball : c^Tx \leq k \right\}.$$
Then the following equalities hold:
\begin{eqnarray}
\nonumber \Delta(X) &=& \dist(\partial \conv(\sphere \setminus X), 0) \\
\nonumber &=& \min(\dist(\partial \ball, 0), \dist(\partial\{x : c^Tx \leq k\}, 0)) \\
\nonumber &=& \min(\dist(\sphere, 0), \dist(\{x : c^Tx = k\}, 0)) \\
\nonumber &=& \min\left(1, \frac{|k|}{\|c\|}\right).
\end{eqnarray}
\HACKENDPROOF

\begin{cor}\label{lemma:deltaMonotone} 
Consider the spherical caps $\mathcal{C}_{c, k_1}$ and $\mathcal{C}_{c, k_2}$ such that $k_1 \leq k_2$. Then we have: $$\Delta(\mathcal{C}_{c,k_1}) \leq \Delta(\mathcal{C}_{c,k_2}).$$
\end{cor}

\begin{lem}\label{lemma:constructSC}
For any set $X \subset \sphere$, there exist $c$ and $k$ such that $\mathcal{C}_{c,k}$ satisfies: $\mathcal{C}_{c,k} \subset X,$ and $\Delta(\mathcal{C}_{c,k}) = \Delta(X).$
\end{lem}

\HACKPROOF 
Let $X' := \conv(S \setminus X)$.
Since $\dist$ is continuous and the set $\partial X'$ is compact, there exists a point $x^* \in \partial X'$, such that:
\begin{eqnarray}\nonumber \Delta(X) = \dist(\partial X', 0) = 
\label{deltaSupporting} \min_{x \in \partial X'}\dist(x, 0) = \dist(x^*, 0).\end{eqnarray} 
Next, consider the hyperplane which is tangent to the ball $\Delta(X)\ball$ at $x^*$, which we denote by $\left\{x : c^Tx = k\right\}$. Then we have:
\begin{equation*}\Delta(X) =  \dist(x^*, 0) = \dist(\{x : c^Tx = k\}, 0) = \min\left(1, \frac{|k|}{\|c\|}\right).
\end{equation*}
Now, consider the spherical cap defined by this tangent plane i.e., $\mathcal{C}_{c, k}$. Then, by Lemma \ref{lemma:delta2} we have
$\Delta(\mathcal{C}_{c,k}) =  \min\left(1, \frac{|k|}{\|c\|}\right)$. Therefore, $\Delta(X) = \Delta(\mathcal{C}_{c,k})$.

\noindent
We next show $\mathcal{C}_{c, k} \subset X$. We prove this by contradiction. Assume $x \in \mathcal{C}_{c,k}$ and $x \notin X$. Note that, if $x \notin X$, then $x \in \sphere \setminus X \subset \conv(\sphere \setminus X).$ Since $x \in \mathcal{C}_{c,k}$, we have $c^Tx>k$. But due to the fact that $x \in \conv(\sphere \setminus X)$, we also have $c^Tx \leq k$, which leads to a contradiction. Therefore, $\mathcal{C}_{c, k} \subset X$. 
\HACKENDPROOF

\noindent
We are now able to prove Proposition 1 given in Section 4 of our paper, which states that, for any $\varepsilon \in (0,1)$, the function $\Delta(X)$ attains its minimum over $\mathcal{X}_{\varepsilon}$ for some $X$ which is a spherical cap, i.e., the minimal radius $\delta$ of the largest sphere $\delta \sphere$ included in $\sphere \setminus X$ will be reached when $X$ is a spherical cap.
\begin{prop}\label{thm:mainSphericalCap}
Let $\tilde{\varepsilon} \in [0,1]$ and $\mathcal{X}_{\tilde{\varepsilon}} = \{X \subset \sphere: \sigma^{n-1}(X) \leq \tilde{\varepsilon}\}$. Then, the function $\Delta(X)$ attains its minimum over $\mathcal{X}_{\tilde{\varepsilon}}$ for some $X$ which is a spherical cap. We denote by $\delta(\tilde{\varepsilon})$ this minimal value, which takes the following expression: $$\delta(\tilde{\varepsilon}) = \sqrt{1 - I^{-1}(2 \tilde{\varepsilon};\frac{n-1}{2},\frac{1}{2})}.$$
\end{prop}

\HACKPROOF
We prove the first part of the proposition via contradiction. Assume that there exists no spherical cap in $\mathcal{X}_{\varepsilon}$ such that $\Delta(X)$ attains its minimum. This means there exists an $X^* \in \mathcal{X}_{\varepsilon}$, where $X^*$ is not a spherical cap and $\argmin_{X \in \mathcal{X}_{\varepsilon}}(\Delta(X))=X^*$. By Lemma \ref{lemma:constructSC}, we can construct a spherical cap $\mathcal{C}_{c,k}$ such that $\mathcal{C}_{c,k} \subset X^*$ and $\mathcal{C}_{c,k} = \Delta(X^*)$. Note that, we further have $\mathcal{C}_{c,k} \subsetneq X^*$, since $X^*$ is assumed not to be a spherical cap. This means that, there exists a spherical cap $\sigma^{n-1}(\mathcal{C}_{c,k})$ such that $\sigma^{n-1}(\mathcal{C}_{c,k}) < \varepsilon$. 

\vspace{0.3cm}

\noindent
Then, the spherical cap $\mathcal{C}_{c, \tilde{k}}$ with $\sigma^{n-1}(\mathcal{C}_{c, \tilde{k}}) = \varepsilon$, satisfies $\tilde{k} < k$ by Remark \ref{lemma:muMonotone}. This implies $$\Delta(\mathcal{C}_{c, \tilde{k}}) < \Delta(\mathcal{C}_{c, k}) = \Delta(X^*)$$ by Corollary \ref{lemma:deltaMonotone}. Therefore, $\Delta(\mathcal{C}_{c, \tilde{k}}) < \Delta(X^*)$. This is a contradiction since we initially assumed that $\Delta(X)$ attains its minimum over $\mathcal{X}_{\varepsilon}$ at $X^*$.

\vspace{0.3cm}

\noindent
We can now give an expression for $\delta(\varepsilon)$. We know that:
\begin{equation}\label{eqn:sc}
\delta(\varepsilon) = \Delta(\mathcal{C}_{c, k}),
\end{equation}
for some spherical cap $\mathcal{C}_{c,k} \subset \sphere$, where  $\sigma^{n-1}(\mathcal{C}_{c, k}) = \varepsilon$. It is known (see e.g. \cite{sphericalCapRef}) that the area of such $\mathcal{C}_{c, k}$, is given by the equation:
\begin{equation}\sigma^{n-1}(\mathcal{C}_{c, k}) = \frac{I\left(1-\Delta(\mathcal{C}_{c,k})^2; \frac{n-1}{2}, \frac{1}{2}\right)}{2}
\end{equation}
where $I$ is the regularized incomplete beta function. Since, \mbox{$\sigma^{n-1}(\mathcal{C}_{c, k})= \varepsilon$,} we get the following set of equations:
\begin{eqnarray}\nonumber \varepsilon &=& \frac{I\left(1- \Delta(\mathcal{C}_{c,k})^2;\frac{n-1}{2}, \frac{1}{2}\right)}{2} \\
\nonumber 1- \Delta(\mathcal{C}_{c, k})^2 &=&  I^{-1}\left(2\varepsilon; \frac{n-1}{2}, \frac{1}{2}\right) \\
\label{eqn:last}\Delta(\mathcal{C}_{c, k})^2 &=&  1- I^{-1}\left(2\varepsilon; \frac{n-1}{2}, \frac{1}{2}\right).
\end{eqnarray}

\noindent
This gives us 
\begin{equation}
\delta(\varepsilon) = \sqrt{1- I^{-1} \left( 2\varepsilon; \frac{n-1}{2}, \frac{1}{2} \right)}.
\end{equation}
\HACKENDPROOF

\bibliographystyle{plain} 
\bibliography{arxiv_version}

\begin{thebibliography}{10}

\bibitem{macLayer}
{I}{E}{E}{E} {Standard} for {Information} {Technology} - {Telecommunications}
  and {Information} {Exchange} {Between} {Systems} - {Local} and {Metropolitan}
  {Area} {Networks} - {Specific} {Requirements} {Part} 15.4: {Wireless}
  {Medium} {Access} {Control} ({M}{A}{C}) and {Physical} {Layer} ({P}{H}{Y})
  {Specifications} for {Low}-{Rate} {Wireless} {Personal} {Area} {Networks}
  ({W}{P}{A}{N}s).
\newblock Technical report, 2006.

\bibitem{balkan}
A.~Balkan, P.~Tabuada, J.~V. Deshmukh, X.~Jin, and J.~Kapinski.
\newblock Underminer: {A} {Framework} for {Automatically} {Identifying}
  {Non}-{Converging} {Behaviors} in {Black} {Box} {System} {Models}.
\newblock In {\em Proceedings of the Thirteenth ACM International Conference on
  Embedded Software (EMSOFT)}, EMSOFT '16, pages 1--10, 2016.

\bibitem{bianchini}
F.~Blanchini, G.~Fenu, G.~Giordano, and F.~A. Pellegrino.
\newblock Model-{Free} {Plant} {Tuning}.
\newblock {\em I{E}{E}{E} {Transactions} on {Automatic} {Control}},
  62(6):2623--2634, 2017.

\bibitem{stabilityHard1}
V.~D. Blondel and J.~N. Tsitsiklis.
\newblock Complexity of {Stability} and {Controllability} of {Elementary}
  {Hybrid} {Systems}.
\newblock {\em Automatica}, 35(3):479--489, 1999.

\bibitem{lazar}
R.~Bobiti and M.~Lazar.
\newblock A {Delta}-sampling {Verification} {Theorem} for {Discrete}-time,
  {Possibly} {Discontinuous} {Systems}.
\newblock In {\em Proceedings of the 18th International Conference on Hybrid
  Systems: Computation and Control}, pages 140--148. ACM, 2015.

\bibitem{boyd}
S.~Boyd and L.~Vandenberghe.
\newblock {\em Convex {Optimization}}.
\newblock Cambridge University Press, 2004.

\bibitem{campi}
G.~Calafiore.
\newblock Random {Convex} {Programs}.
\newblock {\em SIAM Journal on Optimization}, 20(6):3427--3464, 2010.

\bibitem{tiebreak}
G.~Calafiore and M.~C. Campi.
\newblock The {Scenario} {Approach} to {Robust} {Control} {Design}.
\newblock {\em IEEE Transactions on Automatic Control}, 51(5):742--753, 2006.

\bibitem{campi-garatti}
M.~C. Campi and S.~Garatti.
\newblock The {Exact} {Feasibility} of {Randomized} {Solutions} of {Uncertain}
  {Convex} {Programs}.
\newblock {\em SIAM Journal on Optimization}, 19(3):1211--1230, 2008.

\bibitem{wirelessHart}
D.~Chen, M.~Nixon, and A.~Mok.
\newblock {\em WirelessHART: {Real}-{Time} {Mesh} {Network} for {Industrial}
  {Automation}}.
\newblock Springer Publishing Company, Incorporated, 1st edition, 2010.

\bibitem{mitra2}
P.~S. Duggirala, S.~Mitra, and M.~Viswanathan.
\newblock Verification of {Annotated} {Models} from {Executions}.
\newblock In {\em Proceedings of the Eleventh ACM International Conference on
  Embedded Software}, EMSOFT '13, pages 26:1--26:10, 2013.

\bibitem{mitra}
Z.~Huang and S.~Mitra.
\newblock Proofs from {Simulations} and {Modular} {Annotations}.
\newblock In {\em Proceedings of the 17th International Conference on Hybrid
  Systems: Computation and Control}, HSCC '14, pages 183--192. ACM, 2014.

\bibitem{jungers_lncis}
R.~Jungers.
\newblock The {Joint} {Spectral} {Radius}: {Theory} and {Applications}.
\newblock In {\em Lecture Notes in Control and Information Sciences}, volume
  385. Springer-Verlag, 2009.

\bibitem{kapinski}
J.~Kapinski, J.~V. Deshmukh, S.~Sankaranarayanan, and N.~Arechiga.
\newblock Simulation-{Guided} {Lyapunov} {Analysis} for {Hybrid} {Dynamical}
  {systems}.
\newblock In {\em Proceedings of the 17th International Conference on Hybrid
  Systems: Computation and Control}, pages 133--142. ACM, 2014.

\bibitem{kozarev2016case}
A.~Kozarev, J.~Quindlen, J.~How, and U.~Topcu.
\newblock Case {Studies} in {Data}-{Driven} {Verification} of {Dynamical}
  {Systems}.
\newblock In {\em Proceedings of the 19th International Conference on Hybrid
  Systems: Computation and Control}, HSCC '16, pages 81--86. ACM, 2016.

\bibitem{lauer}
F.~Lauer.
\newblock On the {Complexity} of {Switching} {Linear} {Regression}.
\newblock {\em Automatica}, 74:80--83, 2016.

\bibitem{sphericalCapRef}
S.~Li.
\newblock {Concise} {Formulas} for the {Area} and {Volume} of a
  {Hyperspherical} {Cap}.
\newblock {\em Asian Journal of Mathematics and Statistics}, 4:66--70, 2011.

\bibitem{liberzon}
S.~Liu, D.~Liberzon, and V.~Zharnitsky.
\newblock On {Almost} {Lyapunov} {Functions} for {Non}-{Vanishing} {Vector}
  {Fields}.
\newblock In {\em Proceedings of the 55th IEEE Conference on Decision and
  Control (CDC)}, pages 5557--5562, 2016.

\bibitem{betafct}
K.~L. Majumder and G.~P. Bhattacharjee.
\newblock Algorithm {A}{S} 64: {Inverse} of the {Incomplete} {Beta} {Function}
  {Ratio}.
\newblock {\em Journal of the {Royal} {Statistical} {Society}. {Series} {C}
  ({Applied} {Statistics})}, 22:411--414, 1973.

\bibitem{nemirovski}
A.~Nemirovski and A.~Shapiro.
\newblock Convex {Approximations} of {Chance} {Constrained} {Programs}.
\newblock {\em SIAM Journal on Optimization}, 17(4):969--996, 2006.

\bibitem{handbook}
F.~W. Olver.
\newblock {\em NIST {Handbook} of {Mathematical} {Functions} {Hardback} and
  {C}{D}-{R}{O}{M}}.
\newblock Cambridge university press, 2010.

\bibitem{sosLyap}
A.~Papachristodoulou and S.~Prajna.
\newblock On the {Construction} of {Lyapunov} {Functions} {Using} the {Sum} of
  {Squares} {Decomposition}.
\newblock In {\em Proceedings of 41st IEEE Conference on Decision and Control
  (CDC)}, pages 3482--3487, 2002.

\bibitem{parrilo-jadbabaie}
P.~A. Parrilo and A.~Jadbabaie.
\newblock Approximation of the {Joint} {Spectral} {Radius} {Using} {Sum} of
  {Squares}.
\newblock {\em Linear Algebra and its Applications}, 428(10):2385--2402, 2008.

\bibitem{rudin}
W.~Rudin.
\newblock {\em Real and {Complex} {Analysis}}.
\newblock McGraw-Hill Book Co., 3rd edition, 1987.

\bibitem{topcu}
U.~Topcu, A.~Packard, and P.~Seiler.
\newblock Local {Stability} {Analysis} {Using} {Simulations} and
  {Sum}-of-{Squares} {Programming}.
\newblock {\em Automatica}, 44(10):2669--2675, 2008.

\bibitem{jsrtoolbox}
G.~Vankeerberghen, J.~Hendrickx, and R.~Jungers.
\newblock Jsr: {A} {Toolbox} to {Compute} the {Joint} {Spectral} {Radius}.
\newblock In {\em Proceedings of the 17th International Conference on Hybrid
  Systems: Computation and Control}, HSCC '14, pages 151--156. ACM, 2014.

\end{thebibliography}

\end{document}